\input amstex
\documentstyle{amsppt}
\magnification=1200
\hoffset=-0.5pc
\nologo
\vsize=57.2truepc
\hsize=38.5truepc

\spaceskip=.5em plus.25em minus.20em
\define\ldot{ . }

\define\lamz{z}

\define\RRR{R}
\define\ETAA{\eta}

\define\adabotwo{1}
\define\armcusgo{2}
\define\rbieltwo{3}
\define\biemitwo{4}
\define\broetomd{5}
\define\bourbalg{6}
\define\chakirus{7}
\define\cushbate{8}
\define\dalbec{9}
\define\diracone{10}
\define\elhalnee{11}
\define\bhallfou{12}
\define\bhallone{13}
\define\helbotwo{14}
\define\poiscoho{15}
 \define\souriau{16}
 \define\kaehler{17}
   \define\lradq{18}
      \define\qr{19}
   \define\varna{20}
\define\bedlepro{21}
\define\holopewe{22}
\define\hurusch{23}
\define\humphone{24}
\define\kemneone{25}
\define\kirwaboo{26}
\define\lempszoe{27}
\define\lermonsj{28}
\define\lunatwo{29}
\define\lunathr{30}
\define\enelson{31}
\define\netto{32}
\define\procschw{33}
\define\richaone{34}
\define\gwschwar{35}
\define\gwschwat{36}
\define\slodoboo{37}
\define\snowone{38}
\define\steinbtw{39}
\define\szoekone{40}
\define\taylothr{41}
\define\vaccarin{42}
\define\weylbook{43}

\define\holo{complex }

\topmatter
\title Stratified K\"ahler structures on  adjoint quotients
\endtitle
\author Johannes Huebschmann\footnote"*"{Support by the
German Research Council (Deutsche Forschungsgemeinschaft) in the
framework of a Mercator visiting professorship is gratefully
acknowledged. \hfill \hfill}
\endauthor
\affil
Universit\'e des Sciences et Technologies
de Lille
\\
U. F. R. de Math\'ematiques
\\
CNRS-UMR 8524
\\
F-59655 VILLENEUVE D'ASCQ C\'edex, France
\\
Johannes.Huebschmann\@math.univ-lille1.fr
\endaffil
\date{October 30, 2006}
\enddate
\abstract{Given a compact Lie group, endowed with a bi-invariant
Riemannian metric, its complexification inherits a K\"ahler
structure having twice the kinetic energy of the metric as its
potential, and K\"ahler reduction with reference to the adjoint
action yields a stratified K\"ahler structure on the resulting
adjoint quotient. Exploiting classical invariant theory, in
particular bisymmetric functions and variants thereof, we explore
the singular Poisson-K\"ahler geometry of this quotient. Among
other things we prove that, for various compact groups, the real
coordinate ring of the adjoint quotient is generated, as a Poisson
algebra, by the real and imaginary parts of the fundamental
characters. We also show that singular K\"ahler quantization of
the geodesic flow on the reduced level yields the irreducible
algebraic characters of the complexified group.}
\endabstract

\address{\smallskip\noindent USTL, UFR de Math\'ematiques, CNRS-UMR 8524
\newline\noindent
59655 VILLENEUVE D'ASCQ C\'edex, France
\newline\noindent
Johannes.Huebschmann\@math.univ-lille1.fr}
\endaddress
\endkeywords

\endtopmatter
\document
\leftheadtext{Johannes Huebschmann}

\beginsection Introduction

Given a smooth manifold $M$ endowed with an action of a compact
Lie group $K$, the action lifts to a hamiltonian action on the
(total space of the) cotangent bundle $\roman T^*M$ in an obvious
fashion. The reduced space at zero momentum is then a stratified
symplectic space. In recent years various attempts have been made
to understand the singular structure of this kind of reduced
space; for example, unless there is a single stratum, the strata
are not cotangent bundles on strata of the orbit space of $M$. In
this paper we will elucidate the singular structure explicitly for
the special case where $M$ is $K$ itself, endowed with the {\it
conjugation\/} action. A choice of bi-invariant Riemannian metric
amounts to fixing  the kinetic energy,
 and we will  show that
quantization of the reduced kinetic energy on the resulting
singular quotient yields the irreducible algebraic characters of
the complexified group $K^{\Bbb C}$.
%
Reduced spaces of this kind arise in mechanics, and the total
space of the cotangent bundle on a compact Lie group with symmetry
coming from conjugation is the building block for certain lattice
gauge theories. Despite the huge literature on reduction of the
cotangent bundle of a Lie group relative to {\it left
translation\/} (or right translation), the {\it conjugation\/}
action has received little attention. Quantization on the
symplectic quotient (reduced space) of a space of the kind $\roman
T^*K$ at zero momentum relative to the conjugation action
 will provide a step towards understanding
quantization of certain constrained systems.

The polar decomposition of the complexification $K^{\Bbb C}$ of
$K$ and a choice of bi-invariant Riemannian metric on $K$ induce a
bi-invariant diffeomorphism between  $K^{\Bbb C}$ and $\roman
T^*K$ in such a way that the symplectic and complex structures
combine to a K\"ahler structure. Then the reduced space at zero
momentum $(\roman T^*K)_0$ may be identified with the complex
algebraic categorical quotient $K^{\Bbb C}\big/\big/K^{\Bbb C}$
and thereby acquires a complex algebraic structure which, when $K$
is simple and simply connected of rank $r$ (say), comes down to
ordinary $r$-dimensional complex affine space. Indeed, in an
obvious manner, the complex algebraic categorical quotient of
$K^{\Bbb C}$ is isomorphic to the orbit space $T^{\Bbb C}\big / W$
of the complexification $T^{\Bbb C}$ of a maximal torus $T$ in $K$
relative to the action of the Weyl group $W$ on $T^{\Bbb C}$. In
the literature, an orbit space of the kind $T^{\Bbb C}\big / W$ is
referred to as an {\it adjoint quotient\/}. When $K$ is simple and
simply connected of rank $r$ (say), in view of an observation of
{\it Steinberg\/}'s \cite\steinbtw, the fundamental characters
$\chi_1,\dots,\chi_r$ of $K^{\Bbb C}$ furnish a map
 from $K^{\Bbb C}$ onto $r$-dimensional complex affine space $\Bbb A^r$
which identifies the adjoint quotient with $\Bbb A^r$. In
particular, the complex coordinate ring $\Bbb C[T^{\Bbb C}\big /
W]$ is the polynomial algebra in the characters of the fundamental
irreducible representations, and this ring is that of
$W$-invariants of the complex coordinate ring $\Bbb C[T^{\Bbb C}]$
of the maximal torus $[T^{\Bbb C}$ of $K^{\Bbb C}$.

However, the reduced space has more structure: We recall that a
{\it complex analytic stratified K\"ahler space\/} in the sense of
\cite\kaehler\ is a stratified symplectic space
$(N,C^{\infty}(N),\{\,\cdot\,,\,\cdot\,\})$ together with a
compatible complex analytic structure on $N$. That is to say: $N$
comes with (i) a stratification, with (ii) a Poisson algebra
$(C^{\infty}(N),\{\,\cdot\,,\,\cdot\,\})$ of continuous functions,
referred to as a {\it stratified symplectic Poisson algebra\/}
which, on each stratum, restricts to an ordinary smooth symplectic
Poisson algebra, and (iii) with a complex analytic structure, and
the two structures being compatible amounts to the following
additional requirements being satisfied: (iv) each stratum is a
complex analytic subspace which is actually a complex manifold;
(v) holomorphic functions, defined on open subsets of $N$, are
restrictions of functions in $C^{\infty}(N,\Bbb C)=
C^{\infty}(N)\otimes \Bbb C$; and (vi) on each stratum, the
symplectic and complex analytic structures combine to a K\"ahler
structure. Given the compact Lie group $K$, the quotient $(\roman
T^*K)_0$ is a complex analytic stratified K\"ahler space, indeed,
even a \lq\lq complex algebraic\rq\rq\ stratified K\"ahler space
in a sense made precise in Section 1 below. Thus, symplectically
or, more precisely, as a stratified symplectic space, the quotient
$(\roman T^*K)_0$ has singularities, even when $K$ is simple and
simply connected so that complex analytically or complex
algebraically $(\roman T^*K)_0$ is just an affine space. Whether
or not $K$ is simply connected, Poisson brackets among the real
and imaginary parts of holomorphic coordinate functions then yield
(continuous) functions which are not necessarily smooth functions
of the coordinate functions; indeed, the singular structure on the
reduced level is reflected in Poisson brackets among (continuous)
functions which are not necessarily smooth. Explicit examples of
such Poisson brackets will be given in Sections 3 and 4 below.
Suffice it to mention at this stage that the real coordinate ring
$\Bbb R[T^{\Bbb C}\big / W]$ of the quotient $T^{\Bbb C}\big / W$
amounts to the algebra $\Bbb R[T^{\Bbb C}]^W$ of $W$-invariants in
the real coordinate ring $\Bbb R[T^{\Bbb C}]$ of the maximal torus
$T^{\Bbb C}$, viewed as a real algebraic manifold. Now the real
coordinate ring $\Bbb R[T^{\Bbb C}\big / W]$ of the quotient
$T^{\Bbb C}\big / W$ contains of course the subalgebra generated
by the real and imaginary parts of the polynomial generators of
$\Bbb C[T^{\Bbb C}\big / W]$ but these do not generate the real
coordinate ring of the quotient. The stratified symplectic Poisson
structure of the quotient is defined on $\Bbb R[T^{\Bbb C}\big /
W]$, and one of our results is the following.

\proclaim{Theorem} For $K=\roman{U}(n),
\roman{SU}(n),
\roman{Sp}(n),
\roman{SO}(2n+1,\Bbb R),
G_{2(-14)}$,
as a Poisson algebra, the real coordinate ring $\Bbb R[T^{\Bbb C}\big /
W]$ of the quotient $T^{\Bbb C}\big / W$ is generated by the real
and imaginary parts of the characters $\chi_1,\dots,\chi_r$ of the
fundamental irreducible representations of $K^{\Bbb C}$. That is
to say: This ring  is generated by the real and imaginary parts of
these characters, together with iterated Poisson brackets in these
functions.
\endproclaim

For $K=\roman U(n)$ or $K=\roman{SU}(n)$, the theorem comes down
to the statement that the algebra $\Bbb C[z_1,\dots, z_n,
\overline z_1,\dots, \overline z_n]^{S_n}$ of bisymmetric
functions, that is, the algebra of $S_n$-invariants in the
variables $z_1,\dots, z_n, \overline z_1,\dots, \overline z_n$,
where the symmetric group $S_n$ permutes the variables $z_1,\dots,
z_n$ and $\overline z_1,\dots, \overline z_n$ separately, is
generated by the elementary symmetric functions
$\sigma_1,\dots,\sigma_n$ in the variables $z_1,\dots, z_n$ and
the elementary symmetric functions $\overline
\sigma_1,\dots,\overline \sigma_n$ in the variables $\overline
z_1,\dots, \overline z_n$, together with iterated Poisson brackets
in these functions. See Corollaries 3.4.4 and 3.4.8 below for more
details. The cases $\roman{Sp}(n)$, $\roman{SO}(2n+1,\Bbb R)$, and
$G_{2(-14)}$ are established in Section 4.

For $K=\roman{Spin}(n,\Bbb R)$ ($n \geq 3$, $n \neq 4$, $n \neq 6$),
$K=\roman{SO}(2n,\Bbb R)$ ($n \geq 2$),
$K=F_{2(-52)}$, the statement of the theorem is not true.
We do not know what happens for
$K= E_{6(-78)},
E_{7(-132)},
E_{8(-248)}$.

For $K = \roman U(n)$, the unitary group,  complex algebraically,
the reduced space $(\roman T^*K)_0\cong\roman {GL}(n,\Bbb C)\big
/\big/\roman {GL}(n,\Bbb C)$ is the space of complex normalized
degree $n$ polynomials in a single variable having non-zero
constant coefficient, and this space is complex algebraically and
hence complex analytically isomorphic to $\Bbb C^{n-1} \times\Bbb
C^*$ in an obvious fashion. Indeed, the quotient map from $K^{\Bbb
C} \cong\roman {GL}(n,\Bbb C)$ to the space of polynomials sends a
matrix in $\roman {GL}(n,\Bbb C)$ to its characteristic
polynomial, and the stratification of the reduced space is given
by the multiplicities of the roots where strata correspond to
partitions of $n$. The stratified symplectic Poisson algebra
$C^{\infty}(\roman T^*K)_0$ contains more functions than the
ordinary smooth functions on $(\roman T^*K)_0=\Bbb C^{n-1}
\times\Bbb C^*$, though, viewed as a smooth real
$(2n)$-dimensional manifold. In Section 3 below shall determine
the stratified symplectic Poisson structure explicitly. For $K =
\roman {SU}(n)$, the quotient $\roman {SL}(n,\Bbb C)\big/\big/
\roman {SL}(n,\Bbb C)$ amounts to the subspace of the quotient
$\roman {GL}(n,\Bbb C)\big/\big/ \roman {GL}(n,\Bbb C)$ for
$\roman U(n)$ which consists of complex normalized degree $n$
polynomials with constant coefficient equal to $1$. Complex
algebraically, this space is plainly just a copy of $\Bbb
C^{n-1}$.

In the final section we shall show that half-form quantization of
the reduced kinetic energy  associated with the bi-invariant
metric on $K$  yields the irreducible algebraic characters of
$K^{\Bbb C}$. In the situation considered there, quantization
unitarily commutes with reduction. Exploiting results in
\cite\richaone, we plan to extend elsewhere the present approach
to orbit spaces of $n$-tuples of elements from $\roman T^*K$. This
is the typical situation in lattice gauge theory.

This paper was rewritten during a stay at the Institute for
Theoretical Physics at the University of Leipzig. This stay was
made possible by the German Research Council (Deutsche
Forschungsgemeinschaft) in the framework of a Mercator visiting
professorship, and I wish to express my gratitude to this
organization. It is a pleasure to acknowledge the stimulus of
conversation with G. Rudolph and M. Schmidt at Leipzig. The paper
is part of a research program aimed at exploring quantization on
classical phase spaces with singularities \cite\poiscoho
--\cite\bedlepro, in particular on classical lattice gauge theory
phase spaces. Details for the special case of a single spatial
plaquette where $K=\roman{SU}(2)$ are worked out in \cite\hurusch.

\medskip\noindent{\bf 1. The adjoint quotient}
\smallskip\noindent
Let $K$ be a compact Lie group, let $\frak k$ be its Lie algebra,
choose an invariant inner product $\cdot \, \colon \frak k \otimes
\frak k @>>> \Bbb R$ on $\frak k$, and endow $K$ with the
corresponding bi-invariant Riemannian metric. Using this metric,
we identify $\frak k$ with its dual $\frak k^*$ and the total
space of the tangent bundle $\roman T K$ with the total space of
the cotangent bundle $\roman T^* K$. The polar decomposition map
assigns $x\cdot\roman{exp}(iY) \in K^{\Bbb C}$ to $(x,Y) \in K
\times \frak k$. Thus the composite
$$
\roman T^* K @>>> K \times \frak k @>>> K^{\Bbb C} \tag1.1
$$
of the inverse of left trivialization with the polar decomposition
map identifies $\roman T^* K$ with $K^{\Bbb C}$ in a $(K\times
K)$-equivariant fashion. Then the induced complex structure on
$\roman T^* K$ combines with the symplectic structure to a
(positive) K\"ahler structure. Indeed, the {\it real analytic\/}
function
$$
\kappa \colon K^{\Bbb C} @>>> \Bbb R, \quad
\kappa(x\cdot\roman{exp}(iY)) =  |Y|^2,\quad (x,Y) \in K \times
\frak k, \tag1.2
$$
on $K^{\Bbb C}$ which is twice the {\it kinetic energy\/}
associated with the Riemannian metric, is a (globally defined)
{\it K\"ahler potential\/}; in other words, the function $\kappa$
is strictly plurisubharmonic and (the negative of the imaginary
part of) its Levi form yields (what corresponds to) the cotangent
bundle symplectic structure, that is, the cotangent bundle
symplectic structure on $\roman T^*K$ is given by $i \partial
\overline\partial \kappa$. An explicit calculation which
establishes this fact may be found in \cite\bhallone\ (but
presumably it is a folk-lore observation). For related questions
see \cite\lempszoe, \cite\szoekone.

The group $K$ acts on itself and hence on the total space $\roman
T^*K$ of its cotangent bundle via conjugation. The $K$-action on
$\roman T^*K\cong \roman T K$ is hamiltonian and preserves the
K\"ahler structure, with  momentum mapping
$$
\mu \colon \roman TK @>>> \frak k, \quad \mu(X_x) =
X_xx^{-1}-x^{-1}X_x, \tag1.3
$$
where $x \in K$, where $X_x \in \roman T_xK$ is a tangent vector
at $x$, and where $X_xx^{-1}\in \frak k$ and $x^{-1}X_x\in \frak
k$ are the results of right and left translation, respectively,
with $x^{-1}$. By Proposition 4.2 of \cite\kaehler, the K\"ahler
quotient $\roman T^* K\big/\big/K$ at zero momentum is a complex
analytic stratified K\"ahler space.

Symplectically, the quotient is the orbit space $\mu^{-1}(0)\big
/K$, and an observation of {\smc Kempf-Ness} \cite\kemneone\ and
{\smc Kirwan} \cite\kirwaboo, cf. \S 4 of \cite\gwschwat, where
the zero locus $\mu^{-1}(0)$ is referred to as a {\it
Kempf-Ness\/} set, entails that the obvious map
$$
\mu^{-1}(0)\big /K @>>> \roman T^* K\big /\big/K^{\Bbb C} \cong
K^{\Bbb C}\big/\big/ K^{\Bbb C}.
 \tag 1.4
$$
from the symplectic quotient (reduced space) to the categorical
quotient induced by the inclusion of $\mu^{-1}(0)$ into $\roman
T^* K$ is a homeomorphism in the ordinary (not Zariski) topology.
Here $K^{\Bbb C}\big/\big/ K^{\Bbb C}$ refers to
 the complex algebraic categorical quotient of
$K^{\Bbb C}$ relative to the $K^{\Bbb C}$-action on itself via
conjugation; see e.~g. \cite\gwschwat\ (\S 3) for details on the
construction of the categorical quotient in the category of
complex algebraic varieties. In view of results of {\smc Luna\/}
\cite\lunatwo, \cite\lunathr, this quotient is the categorical
quotient in the category of analytic varieties as well, see also
\cite\gwschwar\ (Theorem 3.6).

The categorical quotient $ K^{\Bbb C}\big/\big/ K^{\Bbb C}$, in
turn, has a very simple structure: Choose a maximal torus $T$ in
$K$ and let $W$ be the corresponding Weyl group; then $T^{\Bbb C}$
is a maximal torus in $K^{\Bbb C}$, and the (algebraic) {\it
adjoint quotient\/} $\chi \colon K^{\Bbb C} \to T^{\Bbb C}\big /
W$, cf. \cite\humphone\ (3.4) and \cite\slodoboo\ (3.2) for this
terminology, realizes the categorical quotient. Here $T^{\Bbb
C}\big / W$ is the space of $W$-orbits in $T^{\Bbb C}$, and we
will also refer to the orbit space $T^{\Bbb C}\big / W$ as the
{\it adjoint quotient of \/} $K^{\Bbb C}$. In concrete terms, the
map $\chi$ admits the following description: The closure of the
conjugacy class of $x \in K^{\Bbb C}$ contains a unique semisimple
(equivalently: closed) conjugacy class $C_x$ (say), and semisimple
conjugacy classes are parametrized by $T^{\Bbb C}\big / W$; the
image of $x \in K^{\Bbb C}$ under $\chi$ is simply the parameter
value in $T^{\Bbb C}\big / W$ of the semisimple conjugacy class
$C_x$. Since $W$ is a finite group, as a complex  algebraic space,
the quotient $T^{\Bbb C}\big / W$ is simply the space of
$W$-orbits in $T^{\Bbb C}$.

The choice of maximal torus $T$ in $K$ also provides considerable
simplification for the stratified symplectic structure on the
symplectic quotient $(\roman T^*K)_0=\mu^{-1}(0)\big /K$. Indeed,
via the identification (1.1) for $K=T$, the real space which
underlies the (complex algebraic) orbit space $T^{\Bbb C}\big/W$
for the action of the Weyl group $W$ on $T^{\Bbb C}$ amounts
simply to the orbit space $\roman T^* T\big/ W$, with reference to
the induced action of the Weyl group $W$ on $\roman T^* T$, and
the orbit space $\roman T^* T\big/ W$ inherits a stratified
symplectic structure in an obvious fashion:  Strata are the
$W$-orbits, the closures of the strata are affine varieties, and
the requisite stratified symplectic Poisson algebra
$(C^\infty(\roman T^* T\big/W), \{\,\cdot\,,\,\cdot\,\})$ is
simply the algebra $C^\infty(\roman T^* T)^W$ of smooth
$W$-invariant functions on $\roman T^* T$, endowed with the
Poisson bracket $\{\,\cdot\,,\,\cdot\,\}$ coming from the ordinary
symplectic Poisson bracket on $\roman T^* T$. Moreover, the choice
of invariant inner product on $\frak k$ determines an injection
$\roman T^* T \to \roman T^*K $.

\proclaim{Proposition 1.5} The values of the injection $\roman T^*
T \to \roman T^*K $ lie in the zero locus $\mu^{-1}(0)$, and the
injection $\roman T^* T \to \mu^{-1}(0)$ induces an isomorphism
$$
\roman T^* T\big/W @>>> (\roman T^*K)_0 = \mu^{-1}(0)\big /K
$$
of stratified symplectic spaces.
\endproclaim

\demo{Proof} Since $T$ is an abelian Lie group, and since the
action of $K$ on $\roman T^*K$ is by conjugation, the restriction
of the momentum mapping $\mu$ to $\roman T^* T$, cf. (1.3) above,
is zero, that is, the values of the  injection $\roman T^* T \to
\roman T^*K $ lie in the zero locus $\mu^{-1}(0)$. Since the map
(1.4) and the obvious map from $T^{\Bbb C}\big / W$ to $K^{\Bbb
C}\big / \big /K^{\Bbb C}$ are homeomorphisms, the induced map
from $\roman T^* T\big/W$ to $(\roman T^*K)_0 = \mu^{-1}(0)\big
/K$ is a homeomorphism as well. Moreover, under this map, the
orbit type stratifications correspond.

Any smooth $K$-invariant function on $\roman T^*K$ restricts to a
smooth $W$-invariant function on $\roman T^* T$. Consequently the
homeomorphism from $\roman T^* T\big/W$ to $(\roman T^*K)_0 =
\mu^{-1}(0)\big /K$ induces a map
$$
C^{\infty}((\roman T^*K)_0) =C^{\infty}(\roman T^* K)^K\big /I^K
@>>> (C^{\infty}(\roman T^*T))^W,
$$
necessarily injective since the map between the underlying spaces
is a homeomorphism. Thus is remains to show that each smooth
$W$-invariant function on $\roman T^*T$ extends to a smooth
$K$-invariant function on $\roman T^*K$, with reference to the
conjugation action on $\roman T^*K$. However, this follows from
the compactness of $K$: A smooth $W$-invariant function $f$ on
$\roman T^*T$ extends to a smooth function $\widetilde F$ on
$\roman T^*K$ since $\roman T^*T$ is a closed submanifold of
$\roman T^*T$; averaging over $K$ then yields a smooth
$K$-invariant function $F$ on $\roman T^*K$ extending $f$. \qed
\enddemo

Thus the {\it real\/} structure $C^{\infty}((\roman T^*K)_0)$
comes down to the algebra $C^{\infty}(\roman T^* T)^W$ of (real)
smooth functions on $\roman T^* T \cong T^{\Bbb C}$ that are
invariant under the action of the Weyl group $W$.

The quotient $(\roman T^*K)_0\cong T^{\Bbb C}\big/W$ inherits
various interrelated structures, and for later reference we will
now spell them out and introduce appropriate notation: Given a
real affine locally semialgebraic space $N$ (embedded into some
some real affine space), we write its real coordinate ring as
$\Bbb R[N]$ and, accordingly, we write its ring of analytic
functions and that of Whitney smooth functions as $C^{\omega}(N)$
and $C^{\infty}(N)$, respectively. Likewise, given an affine
complex variety $N$, we denote its complex coordinate ring by
$\Bbb C [N]$. By construction, on the adjoint quotient $N=T^{\Bbb
C}\big/W$, the algebra $\Bbb R[N] = \Bbb R[T^{\Bbb C}]^W$ of
$W$-invariant real polynomial functions on $T^{\Bbb C}$ yields a
real affine locally semialgebraic structure, the algebra $\Bbb
C[N] = \Bbb C[T^{\Bbb C}]^W$ of $W$-invariant complex polynomial
functions on $T^{\Bbb C}$ yields a complex affine structure, and
the real Poisson structure is real algebraic in the sense that it
is defined already on $\Bbb R[N]$. Further, these structures
combine to a {\it complex algebraic stratified K\"ahler
structure\/} on $N$, that is, the Poisson structure is already
defined on $\Bbb R[N]$ and the complex structure is given in terms
of the complex affine coordinate ring $\Bbb C[N]$ but, beware, the
algebra $\Bbb C[N]$ is {\it not\/} the complexification of $\Bbb
R[N]$. Indeed, $\Bbb C[N]$ may be identified with a subalgebra of
the complexification $\Bbb R[N]_{\Bbb C}$ of $\Bbb R[N]$ but $\Bbb
R[N]_{\Bbb C}$ is strictly larger than $\Bbb C[N]$; see Section 3
below for concrete examples. Moreover, the algebra $C^{\omega}(N)
= C^{\omega}(T^{\Bbb C})^W$ of $W$-invariant real analytic
functions on $T^{\Bbb C}$ yields a real affine locally
semianalytic structure on $N$, and the three real structures are
related by the obvious embeddings
$$
\Bbb R[N] \subseteq C^{\omega}(N) \subseteq C^{\infty}(N).
$$
In particular, the K\"ahler potential $\kappa$ (twice the kinetic
energy) on $T^{\Bbb C}$ is a real analytic function which is
plainly $W$-invariant and hence descends to a function
$\kappa_{\roman{red}}$ in $C^{\omega}(N)$, the {\it reduced\/}
K\"ahler potential, which is then twice the reduced kinetic
energy. This function is a K\"ahler potential on the adjoint
quotient $N$ in the sense that, restricted to a stratum, it yields
an ordinary K\"ahler potential on that stratum. Notice that
$\kappa_{\roman {red}}$ does not belong to the real coordinate
ring $\Bbb R[N]$ of the adjoint quotient $N$, though, and the
reduced K\"ahler potential is not an ordinary smooth function,
that is, it is neither a real analytic nor a smooth function on
the adjoint quotient, even when this quotient is topologically
just an affine space. A description of the Poisson bracket on the
real coordinate ring $\Bbb R[N]$ of the adjoint quotient $N$ will
be given in the next section.

It is, perhaps, worthwhile pointing out that, on $T^{\Bbb C}$ and,
likewise, on $K^{\Bbb C}$, the K\"ahler structure is algebraic
while the K\"ahler potential is a real analytic function. On the
other hand,
 the  total spaces $\roman T T$ and $\roman T K$
 of the tangent bundles of the maximal
torus $T$ of $K$ and of $K$ itself, respectively, are complex
analytically equivalent to $T^{\Bbb C}$ and $K^{\Bbb
C}$, respectively, under the polar map, the Poisson structures on
$\roman T T\cong\roman T^* T$ and  on $\roman T K\cong\roman T^*
K$ are real algebraic (even though the identification between
$\roman T^* K$ and $K^{\Bbb C}$ is real analytic), the induced
complex structures on $\roman T T$ and on $\roman T K$ are real
analytic and, on $\roman T^* T$ and $\roman T^* K$, the K\"ahler
potentials are real algebraic functions. Thus we could describe
the stratified K\"ahler structure on the symplectic quotient
$(\roman T^*K)_0$ in terms of the K\"ahler structure on $\roman
T^* T$ as well but this would yield a stratified K\"ahler
structure which is just complex analytic and {\it not\/}
algebraic. On the other hand, the two stratified K\"ahler
structures are plainly equivalent in the category of complex
analytic stratified K\"ahler spaces.

\medskip\noindent{\bf 2. The reduced Poisson algebra}
\smallskip\noindent
As before, let $K$ be a compact Lie group, let $T$ be a maximal
torus of $K$, let $n=\roman{dim}(T) = \roman{rank}(K)$, and let
$W$ be the corresponding Weyl group. In view of Proposition 1.5,
as a stratified symplectic space, the symplectic quotient $(\roman
T^*K)_0= \mu^{-1}(0)\big/K$ amounts to the orbit space $\roman
T^*T\big/W$ and, since $W$ acts symplectically, the Poisson
algebra of smooth $W$-invariant functions on $\roman T^*T\cong
T^{\Bbb C}$ yields the {\it reduced Poisson algebra\/}, that is,
the stratified symplectic Poisson algebra on the quotient. The
purpose of the present section is to derive an explicit
description of this reduced Poisson algebra. Actually, this
Poisson algebra is real algebraic, and we will describe it as a
real algebraic object. To this end we realize the $W$-manifold
$T^{\Bbb C}\cong (\Bbb C^*)^n$ as a closed non-singular
$W$-variety in a suitable $W$-representation. We examine first a
special case, in a manner which may look unnecessarily complicated
but which will pave the way towards quantization on the adjoint
quotient in Section 5 below.

\noindent {\smc (2.1) $K= \roman U(1) = S^1$, the circle group.\/}
We identify the Lie algebra $\roman {Lie}(S^1)$ of
the circle group $S^1$ with the real numbers $\Bbb R$ by means of
the embedding $\Bbb R \to \Bbb C$ given by the association $s
\mapsto -is$,  and we endow $S^1$ with the standard
Riemannian metric. The complexification of the circle group  amounts
to a copy of $\Bbb C^*$. We identify the Lie algebra
$\roman {Lie}(\Bbb C^*)$ of $\Bbb C^*$
with a copy of the complex numbers $\Bbb C$, and
we will use the holomorphic coordinate $w = t + is$ on
this Lie algebra.
 The exponential mapping from $\Bbb C$ to $\Bbb C^*$
factors as
$$
\roman{exp} \colon \Bbb C @>>> S^1 \times \Bbb R @>>> \Bbb C^*
$$
where the first arrow is the universal covering projection (which
sends $w = t + is$ to $(\roman e^{is}, t)$) and the second one the
polar map
$$
S^1 \times \Bbb R \cong \roman T S^1 @>>> \Bbb C^*,\quad
(\roman e^{is}, t)
\mapsto \roman e^t \roman e^{is},\, s, t \in \Bbb R. \tag2.1.1
$$
The K\"ahler potential (1.2) on $S^1 \times \Bbb R$, combined with
the universal covering projection from $\Bbb C$ to $S^1 \times
\Bbb R$, yields the K\"ahler potential $\widetilde \kappa$ on
$\Bbb C$ given by
$$
\widetilde \kappa (w) = t^2 = \left(\frac{w + \overline w}2 \right)^2
= \frac 12 w \overline w +\frac 14 w^2 +\frac 14 \overline w^2;
\tag2.1.2
$$
this K\"ahler potential
yields the standard K\"ahler form on $\Bbb C$
and is manifestly
invariant under the group of deck transformations.
Notice that the standard K\"ahler potential on
$\Bbb C$ (given by the assignment to $w\in \Bbb C$
of $\frac {w \overline w}2$) is {\it not\/} invariant under the group
of deck transformations.
It follows that the induced K\"ahler structure on $\roman T S^1$ is the
ordinary flat one, and the universal covering projection from
$\frak k^{\Bbb C}=\Bbb C$ to $\roman T S^1$ is compatible with the
K\"ahler structures where the universal covering space $\Bbb C$
carries the standard structure. To arrive at an explicit formula
on $\roman T S^1 \cong S^1 \times \Bbb R$, write the Maurer-Cartan
form on $S^1$ as $d s$ and let $z = x+iy$ be the holomorphic
coordinate on (the copy of) $\Bbb C$ (downstairs), so that $\Bbb
C^*\cong \roman T S^1$ appears as a subspace of $\Bbb C$; thus, in
terms of the coordinates $t$ and $s$ appearing in (2.1.1),
$x=\roman e^t \cos (s)$ and $y=\roman e^t \sin (s)$. Under
the polar map, the cotangent bundle symplectic structure $dtds$
on $\roman T S^1 \cong \roman T^* S^1$ (strictly speaking this is
the negative of the ordinary cotangent bundle symplectic
structure) and the standard symplectic structure $dx dy$ on $\Bbb
C^*$ (viewed as a subspace of $\Bbb C$) are related by $dx dy =
\roman e^{2t} dtds$ whence the induced symplectic structure on
$\Bbb C^*$ is given by $\frac 1 {r^2}dx dy$ where $r^2 = z
\overline z$. Indeed, on $\Bbb C^*$, since
$$
t = \frac 12 \log (z \overline z)=\frac 12 \log r^2 = \log r,
$$
the K\"ahler potential $\kappa$ (cf. (1.2)) is given by
$$
\kappa(z) = t^2=\log^2(\sqrt {z \overline z}) =\frac 14 \log^2(z
\overline z) \tag2.1.3
$$
where $\log^2(v) = (\log v)^2$, whence
$$
i \partial \overline\partial \kappa= \frac i {2 r^2}dz \wedge
d\overline z =\frac 1 {r^2}dx\wedge dy . \tag2.1.4
$$
Consequently, in terms of the variables $x$ and $y$, the Poisson
bracket on the real coordinate ring $\Bbb R[\Bbb C^*]$
is given by
$$
\{x,y\} = r^2 \tag2.1.5
$$
where $r^2 = x^2 + y^2$ as usual;
this yields a Poisson structure on
the algebra of smooth real functions $C^{\infty}(\Bbb
C^*)$ on $\Bbb C^*$ in the standard fashion. Alternatively, in terms of the
variables $z$ and $\overline z$, the Poisson bracket on the
complexification $\Bbb R[\Bbb C^*]_{\Bbb C}$ of the
real coordinate ring $\Bbb R[\Bbb C^*]$ is given by
$$
\frac i2 \{z,\overline z\} = z \overline z, \tag2.1.6
$$
and this yields a Poisson structure on
the algebra of smooth complex functions $C^{\infty}(\Bbb
C^*,\Bbb C)$ on $\Bbb C^*$.

\noindent {\smc (2.2) The case of a general compact connected Lie
group\/} $K$: Let $T$ be a maximal torus of $K$, let ${n}=\dim
T=\roman{rank}(K)$, and let $W$ be the Weyl group.

\noindent {\smc (2.2.1) The complex structure.\/} As a
K\"ahler manifold, the complex torus $T^{\Bbb C}$ comes down to a
product $(\Bbb C^*)^{n}$ of $n$ copies of $\Bbb C^*$, each copy of
$\Bbb C^*$ being endowed with the symplectic structure
corresponding to the Poisson structure (2.1.5). The $W$-action on
$T^{\Bbb C}$ induces a $W$-module structure on the complex affine
coordinate ring $\Bbb C[T^{\Bbb C}]$ in an obvious fashion. Choose
coordinate functions $v_1,\dots,v_m$ in $\Bbb C[T^{\Bbb C}]$ which
generate $\Bbb C[T^{\Bbb C}]$ as an algebra  and such that the
$\Bbb C$-linear span of $v_1,\dots,v_m$ in $\Bbb C[T^{\Bbb C}]$ is
$W$-invariant. The assignment to a point $q$ of $T^{\Bbb C}$ of
$$
(v_1(q),\dots,v_m(q))\in V=\Bbb C^{m}
$$
yields an embedding
$$
T^{\Bbb C}\cong(\Bbb C^*)^{n} @>>> V=\Bbb C^{m} \tag2.2.1.1
$$
of  $T^{\Bbb C}\cong(\Bbb C^*)^{n}$ into $V=\Bbb C^{m}$ as a
closed non-singular complex subvariety of $V$. By construction,
$V$ is endowed with a $W$-module structure, and the embedding
(2.2.1.1) is $W$-equivariant. Let $z_1, \dots, z_m$ be the obvious
coordinates on $V$.
 The corresponding morphism
$$
\Bbb C[V]=\Bbb C[z_1,\dots, z_m] @>>> \Bbb C[T^{\Bbb C}]
\tag2.2.1.2
$$
of $\Bbb C$-algebras is given by the assignment to $z_j$ of $v_j$
($1 \leq j \leq m$) and realizes the complex affine coordinate
ring $\Bbb C[T^{\Bbb C}]$ as the quotient of $\Bbb C[V]$ given by
suitable relations. Moreover the $\Bbb C$-algebra $\Bbb C[V]$
inherits a $W$-module structure in an obvious fashion, and the
morphism (2.2.1.2) of $\Bbb C$-algebras is $W$-equivariant. By a
theorem of Hilbert,  the algebra of $W$-invariants $\Bbb C [V]^W$
is finitely generated. This algebra is the complex coordinate ring
$\Bbb C[V/W]$ of the quotient $V/W$, viewed as a complex affine
variety. Since $W$ is finite, the induced map
$$
\Bbb C [V]^W @>>>  \Bbb C[T^{\Bbb C}]^W \tag2.2.1.3
$$
to the algebra $\Bbb C[T^{\Bbb C}]^W$ of $W$-invariants in $\Bbb
C[T^{\Bbb C}]$ is surjective. This  algebra of invariants is the
complex coordinate ring $\Bbb C[T^{\Bbb C}/W]$ of the quotient
$T^{\Bbb C}/W$, and the surjection (2.2.1.3) is dual to the
induced embedding of $T^{\Bbb C}/W$ into $V/W$. Thus a choice
$f_1,\dots,f_k$ of multiplicative generators of $\Bbb C [V]^W$
induces embeddings $ T^{\Bbb C}/W \subseteq V/W \subseteq \Bbb C^k
$ which realize $T^{\Bbb C}/W$ and $V/W$ as complex affine
varieties in $\Bbb C^k$.

\noindent {\smc (2.2.2) The real semialgebraic structure.\/} For
$1 \leq j \leq m$, let $z_j = x_j + y_j$, let $\Bbb R [V]=\Bbb R
[x_1,y_1,\dots, x_m,y_m]$, the real coordinate ring of $V$ where
$V$ is viewed as a real $(2m)$-dimensional vector space and,
 $T^{\Bbb C}$ being viewed as a real $(2n)$-dimensional
 non-singular variety, let
$\Bbb R[T^{\Bbb C}] $ denote its real coordinate ring. By
construction, $\Bbb R[T^{\Bbb C}]$ is multiplicatively generated
by $\roman{Re}(v_1),\roman{Im}(v_1), \dots,
\roman{Re}(v_m),\roman{Im}(v_m)$, and  the assignment to $x_j$ and
$y_j$ of $\roman{Re}(v_j)$ and $\roman{Im}(v_j)$, respectively,
($1 \leq j \leq m$) yields a $W$-equivariant surjection from $\Bbb
R [V]$ onto $\Bbb R[T^{\Bbb C}]$ which is dual to the
$W$-equivariant embedding (2.2.1.1) of $T^{\Bbb C}$ into $V$, both
spaces being viewed as real (non-singular affine) $W$-varieties.
By the theorem of Hilbert quoted earlier,  the algebra of real
$W$-invariants $\Bbb R [V]^W$ is finitely generated. This algebra
is the {\it real\/} coordinate ring $\Bbb R [V/W]$ of the quotient
$V/W$, viewed as a real semialgebraic space (see below). Since $W$
is finite, the induced map
$$
\Bbb R [V]^W @>>>  \Bbb R[T^{\Bbb C}]^W \tag2.2.2.1
$$
to the algebra $\Bbb R[T^{\Bbb C}]^W$ of $W$-invariants in $\Bbb
R[T^{\Bbb C}]$ is surjective. This  algebra of invariants is the
{\it real\/} coordinate ring $\Bbb R[T^{\Bbb C}/W]$ of the
quotient $T^{\Bbb C}/W$, and the surjection (2.2.2.1) is dual to
the induced embedding of $T^{\Bbb C}/W$ into $V/W$, both spaces
being viewed as real semialgebraic sets. These quotients are not
ordinary real varieties, that is, neither of them is the space of
real points of a complex variety. To spell out the structure
somewhat more explicitly, pick a finite system of generators
$\alpha_1,\dots, \alpha_{\ell}$ of $\Bbb R [V]^W$ and let
$$
\alpha = (\alpha_1,\dots,  \alpha_{\ell}) \colon V @>>> \Bbb
R^{\ell}.
$$
This map induces an embedding of the quotient $V/W$ into  $\Bbb
R^{\ell}$ which realizes $V/W$
 as a real {\it semialgebraic\/} set in
$\Bbb R^{\ell}$ \cite\biemitwo\ (\S 1). Defining relations for the
algebra $\Bbb R[V]^W$ of $W$-invariants then yield defining
equations for the smallest {\it real\/} variety $\widehat {V/W}$
in which $V/W$ lies, that is, for the corresponding real {\it
categorical\/} quotient $\widehat {V/W}$ of $V$; as a subspace of
this categorical quotient, the quotient $V/W$ is then given by a
finite set of inequalities which, in turn, encapsulate the real
semialgebraic structure of $V/W$. The real categorical quotient
$\widehat {T^{\Bbb C}\big/W}$ of $T^{\Bbb C}$ and the quotient
$T^{\Bbb C}\big/W$ which are our primary objects of interest are
related in the same fashion: The image of $T^{\Bbb C}$ in $V$ can
be described by a single $W$-invariant equation (by suitable sums
of squares if need be), and the categorical quotient $\widehat
{T^{\Bbb C}\big/W}$ embeds into $\widehat {V/W}$ as a real
algebraic subset of  $\widehat {V/W}$. Hence $T^{\Bbb C}/W \cong
({V/W} \cap \widehat {T^{\Bbb C}\big/W} )\subseteq \widehat
{V/W}$. Thus the inequalities determining ${V/W}$ in $\widehat
{V/W}$ determine $T^{\Bbb C}/W$ as a semialgebraic subset of
$\widehat {T^{\Bbb C}\big/W} $. See e.~g. \cite\gwschwat\ for more
details.

\noindent {\smc (2.2.3) The real algebraic Poisson structure.\/}
Identify $T^{\Bbb C}$ with the product of $n$ copies of $\Bbb
C^*$, each copy of $\Bbb C^*$ being endowed with the Poisson
structure (2.1.5); now, calculating the brackets among the images
$\widetilde\alpha_1,\dots, \widetilde\alpha_{\ell} \in\Bbb
R[T^{\Bbb C}]^W$ of the generators $\alpha_1,\dots, \alpha_{\ell}$
of the algebra $\Bbb R[V]^W$ of $W$-invariants, we obtain a
description of the resulting Poisson structure on the real
coordinate ring $\Bbb R[T^{\Bbb C}/W]$.

Alternatively, we may describe the resulting Poisson structure on
$\Bbb R[T^{\Bbb C}/W]$ in terms of the complexification $\Bbb
R[T^{\Bbb C}/W]_{\Bbb C}$. This provides considerable
simplification, as we shall illustrate in Section 3 below. Indeed,
the complexification $\Bbb R [V]_{\Bbb C}$ of $\Bbb R [V]$
(beware: this complexification is not the algebra $\Bbb C [V]$
considered above) is the polynomial algebra
$$
\Bbb R [V]_{\Bbb C} = \Bbb C[z_1,\overline z_1,\dots, z_m,
\overline z_m] \tag2.2.3.1
$$
in the variables $z_1,\overline z_1,\dots, z_m, \overline z_m$, as
a $\Bbb C$-algebra, the complexification $\Bbb R[T^{\Bbb C}]_{\Bbb
C}$ of $\Bbb R[T^{\Bbb C}]$ is multiplicatively generated by
$v_1,\overline v_1, \dots, v_m,\overline v_m$, the algebras $\Bbb
R [V]_{\Bbb C} $ and $\Bbb R[T^{\Bbb C}]_{\Bbb C}$ carry induced
$W$-module structures in an obvious fashion, and the assignment to
$z_j$ of $v_j$ and to $\overline z_j$ of $\overline v_j$ ($1 \leq
j \leq m$) yields the corresponding $W$-equivariant surjection
from $\Bbb R [V]_{\Bbb C} $ onto $\Bbb R[T^{\Bbb C}]_{\Bbb C}$.
Furthermore,  the algebra of $W$-invariants $\Bbb R [V]_{\Bbb
C}^W$ is still finitely generated. This algebra is the
complexification of the {\it real\/} coordinate ring of the
quotient $V/W$, viewed as a real semialgebraic space. Since $W$ is
finite, the induced map
$$
\Bbb R [V]_{\Bbb C}^W = \Bbb C[z_1,\overline z_1,\dots, z_m,
\overline z_m]^W @>>>  \Bbb R[T^{\Bbb C}]_{\Bbb C}^W \tag2.2.3.2
$$
to the algebra $\Bbb R[T^{\Bbb C}]_{\Bbb C}^W$ of $W$-invariants
in $\Bbb R[T^{\Bbb C}]_{\Bbb C}$ is surjective. This  algebra of
invariants is the complexification $\Bbb R[T^{\Bbb C}/W]_{\Bbb C}$
of the real coordinate ring $\Bbb R[T^{\Bbb C}/W]$  of the
quotient $T^{\Bbb C}/W$. Similarly as before, we may now pick a
finite system of generators $\beta_1,\dots, \beta_{\ell}$ of $\Bbb
R [V]^W_{\Bbb C} = \Bbb C[z_1,\overline z_1,\dots, z_m,\overline
z_m]^W$. In particular, once a choice $f_1,\dots,f_k$ of
multiplicative generators of $\Bbb C [V]^W$, referred to
henceforth as {\it \holo\/} invariants, has been made, cf. what
was said above, these \holo generators and their complex
conjugates yield $2k$  invariants in $\Bbb C[z_1,\overline
z_1,\dots, z_m,\overline z_m]^W$; however, these $2k$ invariants
will not generate the algebra of invariants, and the system of
invariants must be completed by {\it mixed\/} invariants,
that is, by invariants involving the $z_j$'s and the $\overline
z_j$'s ($1 \leq j \leq k$).

\noindent {\smc (2.2.4) The smooth Poisson structure.\/} In view
of a result of {\smc G. W. Schwarz\/} \cite\gwschwar, every smooth
$W$-invariant function of the variables $x_1,y_1,\dots, x_m,y_m$
can be written as a smooth function of the variables
$\alpha_1,\dots, \alpha_{\ell}$, that is, under the projection
from $V$ to $V/W$, the algebra of Whitney-smooth functions on
$V/W$ (relative to the embedding into  $\Bbb R^{\ell}$) is
identified with the algebra $C^{\infty}(V)^W$ of smooth
$W$-invariant functions on $V$. Furthermore, since $W$ is finite,
the induced map
$$
C^{\infty}(V)^W@>>> C^{\infty}(T^{\Bbb C})^W
$$
is surjective whence every smooth $W$-invariant function on
$T^{\Bbb C}$ can be written as a smooth $W$-invariant function of
the variables $x_1,y_1,\dots, x_m,y_m$  and hence as a smooth
function in the variables $\alpha_1,\dots, \alpha_{\ell}$.
Consequently, under the projection from $T^{\Bbb C}$ to $T^{\Bbb
C}/W$, the algebra of real Whitney-smooth functions on $T^{\Bbb
C}/W$ (relative to the embedding into  $\Bbb R^{\ell}$) is
identified with the algebra $C^{\infty}(T^{\Bbb C})^W$ of real
smooth $W$-invariant functions on $T^{\Bbb C}$.

The Poisson structure on the real coordinate ring $\Bbb R[T^{\Bbb
C}/W]$ induces the reduced Poisson structure
$\{\,\cdot\,,\,\cdot\,\}$ on the algebra $C^{\infty}(N)$ in terms
of generators of the algebra $C^{\infty}(T^{\Bbb C}/W) \cong
C^{\infty}(T^{\Bbb C})^W$ as promised at the beginning of the
present section. Likewise, in view of the quoted result of {\smc
G. W. Schwarz\/}, every smooth complex $W$-invariant function of the
variables $z_1,\overline z_1,\dots, z_m, \overline z_m$ can be
written as a smooth function in the variables $\beta_1,\dots,
\beta_{\ell}$, and, likewise, every smooth complex $W$-invariant
function on $T^{\Bbb C}$ can be written as a smooth $W$-invariant
function of the variables $z_1,\overline z_1,\dots, z_m, \overline
z_m$  and hence as a smooth function in the variables
$\beta_1,\dots, \beta_{\ell}$. Now, with reference to the product
Poisson bracket on $C^{\infty}(T^{\Bbb C},\Bbb C)$ where $T^{\Bbb
C}$ is identified with the product of $n$ copies of $\Bbb C^*$,
each copy of $\Bbb C^*$ being endowed with the complex Poisson
structure (2.1.6), calculating the brackets among the images
$\widetilde\beta_1,\dots, \widetilde\beta_{\ell} \in\Bbb R[T^{\Bbb
C}]_{\Bbb C}^W$ of the generators $\beta_1,\dots, \beta_{\ell}$ of
the algebra $\Bbb R[V]_{\Bbb C}^W$ of $W$-invariants, we obtain a
description of the reduced Poisson algebra $(C^{\infty}(T^{\Bbb C}\big/W,\Bbb
C),\{\,\cdot\,,\,\cdot\,\})$ in terms of generators of the algebra
$C^{\infty}(T^{\Bbb C}\big/W,\Bbb C) \cong C^{\infty}(T^{\Bbb C},\Bbb C)^W$.
Since $T^{\Bbb C}$ is actually a K\"ahler manifold, once a choice
$f_1,\dots,f_k$ of \holo invariants has been made, Poisson
brackets of the kind $\{f_u, f_v\}$ and $\{\overline f_u,
\overline f_v\}$ ($1 \leq u,v \leq k$) are necessarily zero.

\medskip\noindent{\bf 3. Unitary and special unitary groups}
\smallskip\noindent
We begin with the unitary group $K= \roman U(n)$.

\noindent {\smc (3.1) The complex structure of the adjoint
quotient.\/} The complexification $K^{\Bbb C}$ of $K= \roman U(n)$
is the full general linear group  $\roman {GL}(n,\Bbb C)$, the
standard maximal torus $T^{\Bbb C} \cong (\Bbb C^*)^n$ consists of
the diagonal matrices in $\roman {GL}(n,\Bbb C)$, and the Weyl
group $W$ is the symmetric group $S_n$ on $n$ letters which acts
on $(\Bbb C^*)^n$ by permutation of the factors. Since the action
of the Weyl group on $T^{\Bbb C}\cong (\Bbb C^*)^n$ actually
extends to the ambient complex vector space $\Bbb C^n$, the
complex affine coordinate ring $\Bbb C[T^{\Bbb C}]$ may be written
as
$$
\Bbb C[T^{\Bbb C}] = \Bbb C[v_1,\dots, v_n,\sigma_n^{-1}]
$$
where $\sigma_n = v_1\cdot \dots \cdot v_n$ is the $n$'th
elementary symmetric function, in such a way that the $W$-action
permutes the $v_1,\dots, v_n$ and leaves invariant the coordinate
function $\sigma_n^{-1}$. Thus we may carry out the construction
in (2.2) above with $m=n+1$, $v_m = \sigma_n^{-1}$, and $V= \Bbb
C^{n+1}$ the $W$-representation  $V=\Bbb C^n \oplus \Bbb C$ where
$W$ acts on the copy of $\Bbb C^n$ by permutation of the factors
and trivially on the residual copy $\Bbb C$. This yields $T^{\Bbb
C}$ as the complex $W$-subvariety of $\Bbb C^{n+1}$ given by the
single equation
$$
v_1\cdot \dots \cdot v_n \cdot v_{n+1} = 1.
$$
This observation, in turn, implies at once that the adjoint
quotient may be realized by means of the map
$$
(\sigma_1,\dots,\sigma_n) \colon (\Bbb C^*)^n @>>> \Bbb C^{n-1}
\times \Bbb C^* \tag3.1.1
$$
whose target identifies the adjoint quotient $T^{\Bbb C}\big/S_n$
as a complex algebraic and hence complex analytic space (which may
be realized as a complex affine variety in $\Bbb C^{n+1}$ as
explained above). Here the constituents $\sigma_1,\dots,\sigma_n$
are the elementary symmetric functions; in fact, these are the
restrictions to $T^{\Bbb C}$ of the characters of the fundamental
(finite dimensional) representations of $\roman{GL}(n,\Bbb C)$.

After incorporation of the appropriate signs, which amounts to
multiplying the $j$'th elementary symmetric function by $(-1)^j$,
the description of the quotient map from $\roman {GL}(n,\Bbb C)$
to the adjoint quotient in terms of polynomials given in the
introduction follows immediately from these observations. The
stratum corresponding to the partition $n_1 + \dots + n_k = n$
consists of the polynomials
$$
(\lamz - z_1)^{n_1} \dots (\lamz - z_k)^{n_k}
$$
in the variable $z$
where $z_j\ne z_{\ell}$ when $j \ne \ell$ ($1 \leq j,\ell \leq k$);
this stratum is a $k$-dimensional complex algebraic manifold. Indeed, the
symmetric group $S_k$ on $k$ letters
acts on the space of all
$$
(z_1, \dots, z_1, z_2, \dots, z_2, \dots, z_k, \dots, z_k),\quad
z_j\ne z_{\ell}, j \ne \ell \ (1 \leq j,\ell \leq k),
$$
in an obvious manner, and the projection to the $S_n$-quotient
identifies the $S_k$-quotient with the corresponding stratum. In
particular, the top stratum, that is, the open and dense stratum
of complex dimension $n$, corresponds to the partition $n_1 +
\dots + n_k = n$ of $n$ where $k=n$ and where each $n_j=1$ ($1
\leq j \leq n$). The complement of the top stratum is, then, the
{\it discriminant\/} variety.  This is the variety over which the
projection mapping (3.1.1) branches; it is given by the equation
$$
D_n(1, -\sigma_1,\sigma_2, \dots, (-1)^n \sigma_n) = 0
$$
where $D_n(a_0,a_1,\dots,a_n)$ refers to the discriminant of the polynomial
$$
P(w) = a_0 w^n + a_1 w^{n-1} + \dots + a_n.
$$
A detailed description of the closures of the strata as
discriminant varieties may be found in \cite\bedlepro.

\noindent {\smc (3.2) The real coordinate ring of the adjoint
quotient.\/} We will spell out a description of the {\it
complexification\/} $\Bbb R [T^{\Bbb C}\big/S_n]_{\Bbb C}$ of the
real coordinate ring $\Bbb R [T^{\Bbb C}\big/S_n]$ of the adjoint
quotient  $T^{\Bbb C}\big/S_n$ of $K^{\Bbb C}=\roman {GL}(n,\Bbb
C)$, where $T^{\Bbb C}$ is the complexification of the maximal
torus $T$ in $K=\roman U(n)$.

Under the present circumstances, the complex algebra $\Bbb R
[V]_{\Bbb C}$, cf. (2.2.3.1), is the polynomial algebra
$$
\Bbb R [V]_{\Bbb C} = \Bbb C[z_1,\overline z_1, \dots, z_{n+1},
\overline z_{n+1}]
$$
in the variables $z_1,\overline z_1,\dots, z_{n+1}, \overline
z_{n+1}$, and the $S_n$-action on $\Bbb R [V]_{\Bbb C}$ permutes
the $z_1,\dots, z_n$'s and the $\overline z_1,\dots, \overline
z_n$'s separately and fixes $z_{n+1}$ and $\overline z_{n+1}$.
Hence the algebra of invariants $\Bbb R [V]_{\Bbb C}^{S_n}$ may be
written as
$$
\Bbb R [V]_{\Bbb C}^{S_n} \cong \Bbb C[z_1,\overline z_1,\dots,
z_n, \overline z_n]^{S_n} \otimes \Bbb C[z_{n+1}, \overline
z_{n+1}].
$$
Thus it suffices to determine the algebra $\Bbb C[z_1,\overline
z_1,\dots, z_n,\overline z_n]^{S_n}$ of $S_n$-invariants which, in
turn, is the algebra of invariants on the product of two copies of
the standard permutation representation of the symmetric group
$S_n$. We will refer to an algebra of this kind as an {\it algebra
of bisymmetric functions\/}.

Systems of generators for an algebra of bisymmetric functions are
well known and classical. To recall such a system, we will follow
Ex. 5 in \S 5 of Chap. IV in \cite\bourbalg. Consider the
polynomial ring $\Bbb Q[ z_1,\dots, z_n,\overline z_1, \dots,
\overline z_n]$ over the field $\Bbb Q$ of rational numbers,
endowed with the obvious $S_n$-action which permutes the variables
separately. For $r \geq 0$ and $s \geq 0$ such that $1 \leq r + s
\leq n$, let
$$
\sigma_{(r,s)}(z_1,\dots, z_n,\overline z_1, \dots, \overline z_n)
\tag3.2.1
$$
be the $S_n$-orbit sum of the monomial $ z_1\cdot\ldots\cdot z_r
\overline z_{1}\cdot\ldots\cdot \overline z_{s}$; the function
$\sigma_{(r,s)}$ is manifestly $S_n$-invariant, and we refer to a
function of this kind as an {\it elementary bisymmetric
function\/}. In degree $m$, $1 \leq m \leq n$, the construction
yields the $m+1$ bisymmetric functions $\sigma_{(m,0)}$,\
$\sigma_{(m-1,1)}$, \dots, $\sigma_{(0,m)}$ whence in degrees at
most equal to $n$  it yields altogether $\frac {n(n+3)}2$
elementary bisymmetric functions; in particular,
$$
\aligned
\sigma_{(m,0)}
(z_1,\dots, z_n,\overline z_1, \dots, \overline z_n) &=
\sigma_m(z_1,\dots, z_n),\\
\sigma_{(0,m)}
(z_1,\dots, z_n,\overline z_1, \dots, \overline z_n) &=
\sigma_m(\overline z_1, \dots, \overline z_n).
\endaligned
$$
According to a classical result, the ring $\Bbb Q[ z_1,\dots,
z_n,\overline z_1, \dots, \overline z_n]^{S_n}$ of bisymmetric
functions (over the rationals) is generated by the elementary
bisymmetric functions \cite\bourbalg. The result is actually more
general and refers to general multisymmetric functions. Suitable
multiples of the elementary bisymmetric functions arise from the
elementary symmetric functions by {\it polarization\/} and, in
terms of the polarizations of the elementary symmetric functions,
the result, for general multisymmetric functions, is given  in
II.3 (p.~37) of \cite\weylbook. Modern accounts may be found in
\cite\dalbec,\ \cite\elhalnee, \ \cite\vaccarin.

The generators (3.2.1) will provide a satisfactory description of
the real semialgebraic geometry of the adjoint quotient. On the
other hand, we shall see in Theorem 3.4.1 below that the
stratified symplectic Poisson algebra can much more easily be
described in terms of the following system of multiplicative
generators for the algebra of bisymmetric functions, which is
entirely classical as well: For $r \geq 0$ and $s \geq 0$ such
that $1 \leq r + s$, let
$$
\tau_{(r,s)}(z_1,\dots, z_n,\overline z_1, \dots, \overline z_n) =
\sum_{j=1}^n z^r_j\overline z^s_j; \tag3.2.2
$$
plainly, $\tau_{(r,s)}$ is the $S_n$-orbit sum of the monomial $
z^r_1\overline z^s_1$, such a function is manifestly
$S_n$-invariant, and we refer to a function of this kind as a {\it
bisymmetric power sum function\/}. Then $\tau_{(1,0)} =
\sigma_{(1,0)}=\sigma_1$ and $\tau_{(0,1)} =
\sigma_{(0,1)}=\overline \sigma_1$ and, in degree $m$, $2 \leq m
\leq n$, the construction yields the $m+1$ bisymmetric functions
$\tau_{(m,0)}$,\ $\tau_{(m-1,1)}$,\ \dots, $\tau_{(0,m)}$ whence
it yields altogether $\frac {n(n+3)}2$ bisymmetric power sum
functions; in particular, when $\tau_m(w_1,\dots, w_n)= w_1^m +
\ldots + w_n^m$ denotes the ordinary $m$'th power sum function,
$$
\aligned \tau_{(m,0)} (z_1,\dots, z_n,\overline z_1, \dots,
\overline z_n) &= \tau_m(z_1,\dots, z_n)= z_1^m + \ldots + z_n^m
,\\
\tau_{(0,m)} (z_1,\dots, z_n,\overline z_1, \dots, \overline z_n)
&= \tau_m(\overline z_1, \dots, \overline z_n)= \overline z_1^m +
\ldots + \overline z_n^m,
\endaligned
$$
and we will occasionally write $\tau_m$ instead of $\tau_{(m,0)}$
and $\overline \tau_m$ instead of $\tau_{(0,m)}$.

According to a variant of the already quoted classical result,
over the rationals, the ring of bisymmetric functions in the
variables $z_1$, \dots, $z_n$, $ \overline z_1$, \dots, $\overline
z_n$  is generated by the bisymmetric power sum functions
$\tau_{(r,s)}$ with $1 \leq r+s \leq n$ as well \cite\dalbec,\
\cite\elhalnee, \ \cite\vaccarin.

Neither the generators (3.2.1) nor the generators (3.2.2) are
free, that is, the algebra of bisymmetric functions is not a
polynomial algebra when $n>1$. It seems hard to find defining
relations in the literature, in fact, we do not know of any
reference for such relations except the description in Theorem 2
(1) of \cite\vaccarin\ where a procedure for constructing
relations in terms of an appropropriate infinite system of
generators is given. We will therefore explain how a system of
defining relations can be derived from a classical result that has
been known since the 19'th century, cf. \cite\netto.

Maintaining standard terminology, we will refer to the quotient
field of the ring of bisymmetric functions (over the rationals)
 as  the {\it field of
bisymmetric functions\/} (over the rationals) . This field is
plainly the field $\Bbb Q(z_1,\dots,z_n,\overline
z_1,\dots,\overline z_n)^{S_n}$ of invariants in the field $\Bbb
Q(z_1,\dots,z_n,\overline z_1,\dots,\overline z_n)$ of rational
functions in $z_1$,\dots,$z_n$, $\overline z_1$,\dots,$\overline
z_n$. The following result is classical, cf. \cite\netto.

\proclaim{Proposition 3.2.3} The field $\Bbb
Q(z_1,\dots,z_n,\overline z_1,\dots,\overline z_n)^{S_n}$ of
bisymmetric functions is purely transcendantal of degree $2n$,
having the functions
$$
\sigma_{(1,0)},\,\dots,\sigma_{(n,0)},\,\sigma_{(0,1)},\,
\sigma_{(1,1)},\,\dots,\sigma_{(n-1,1)} \tag3.2.4
$$
as free generators. In particular, each bisymmetric function can
be written in a unique fashion as a rational function with
rational coefficients in the functions {\rm (3.2.4)}. Likewise,
the functions
$$
\tau_{(1,0)},\,\dots,\tau_{(n,0)},\,\tau_{(0,1)},\,
\tau_{(1,1)},\,\dots,\tau_{(n-1,1)} \tag3.2.5
$$
are free generators of the field of bisymmetric functions, whence
 each bisymmetric function can be written in a unique fashion as
a rational function with rational coefficients in the functions
{\rm (3.2.5)}.
\endproclaim

This is a special case of a more general classical result
saying that the field of $m$-multi-$S_n$-symmetric functions is purely
transcendental with $mn$ generators; cf. \cite\netto. This result,
in turn, is a special case of the familiar fact that any symmetric power
of a rational variety is again rational.

We sketch a proof, since we shall refer to
the argument later in the paper.

\demo{Proof} We will only justify the claim involving the
elementary bisymmetric functions.

In terms of the functions
$$
\sigma_{(0,1)},\,
\sigma_{(1,1)},\,\dots,\sigma_{(n-1,1)},z_1,\dots, z_n, \tag3.2.6
$$
for $1 \leq j \leq n$, the generator $\overline z_j$  is given by
the expression
$$
\overline z_j = \frac {z^{n-1}_j \sigma_{(0,1)} -
z^{n-2}_j\sigma_{(1,1)} + \ldots + (-1)^{n-1}\sigma_{(n-1,1)}}
{(z_j-z_1)\cdot\ldots\cdot (z_j-z_n)} \tag3.2.7
$$
where the denominator is a product of $n-1$ non-zero terms, that
is, the $j$'th term which would formally be of the kind $z_j-z_j$
does {\it not\/} occur. This yields the $\overline z_j$'s ($1 \leq
j \leq n$) as rational functions of the functions (3.2.6).
Consequently the subfield of the field $\Bbb
Q(z_1,\dots,z_n,\overline z_1,\dots,\overline z_n)$ generated by
the functions (3.2.6) exhausts this field whence the field $\Bbb
Q(z_1,\dots,z_n,\overline z_1,\dots,\overline z_n)^{S_n}$ of
bisymmetric functions coincides with the field
$$
\Bbb Q(\sigma_{(0,1)},\,
\sigma_{(1,1)},\,\dots,\sigma_{(n-1,1)},z_1,\dots, z_n)^{S_n}
$$
which, in turn, comes down to the field
$$
\Bbb Q(\sigma_{(0,1)},\,
\sigma_{(1,1)},\,\dots,\sigma_{(n-1,1)},\sigma_1,\dots, \sigma_n).
$$
Since the field $\Bbb Q(z_1,\dots,z_n,\overline
z_1,\dots,\overline z_n)$ has transcendence degree $2n$ and since
it is a Galois extension of the field $\Bbb Q(\sigma_{(0,1)},\,
\sigma_{(1,1)},\,\dots,\sigma_{(n-1,1)},\sigma_1,\dots, \sigma_n)$
of bisymmetric functions, the latter has transcendence degree $2n$
as well, and it is in fact the field of rational functions in the
generators $\sigma_{(0,1)}$, \dots, $\sigma_{(n-1,1)}$,
$\sigma_1$, \dots, $\sigma_n$.

The assertion involving the bisymmetric power sum functions is
established in a similar fashion. \qed \enddemo

Thus, for $2 \leq s \leq n$ and $0 \leq r \leq n$ with $r+s \leq n$,
each generator $\sigma_{(r,s)}$
can be written
in a unique fashion as a rational function
$$
\sigma_{(r,s)} = \frac{\alpha_{(r,s)}}{\beta_{(r,s)}} \tag3.2.8
$$
with rational coefficients, where the $\alpha_{(r,s)}$'s and
$\beta_{(r,s)}$'s are polynomials in the functions {\rm (3.2.4)}
with rational coefficients such that each  $\alpha_{(r,s)}$ is
relatively prime to $\beta_{(r,s)}$. Indeed, each of the rational
functions on the right-hand side of (3.2.8) may be determined in
the following fashion: Given $2 \leq s \leq n$ and $0 \leq r \leq
n$ with $r+s \leq n$, in the definition (3.2.1) of
$\sigma_{(r,s)}$, for $1 \leq j \leq n$, substitute the right-hand
side of (3.2.7) for each occurrence of $\overline z_j$. This
yields each $\sigma_{(r,s)}$ as an $S_n$-invariant rational
function of the functions (3.2.6), that is, when $L$ denotes the
field $\Bbb Q(\sigma_{(0,1)},\,
\sigma_{(1,1)},\,\dots,\sigma_{(n-1,1)})$ of rational functions in
$\sigma_{(0,1)}$,\ $\sigma_{(1,1)}$,\ \dots, $\sigma_{(n-1,1)}$,
the procedure yields $\sigma_{(r,s)}$ as an $S_n$-invariant
function in $L(z_1,\dots,z_n)$, i.~e. as a uniquely determined
rational function in $L(\sigma_1,\dots,\sigma_n)$. In practice,
determining such a rational function explicitly may be a difficult
endeavour, though; see  (3.7) and (3.10) below.

\proclaim{Corollary 3.2.9} For $s \geq 2$, the resulting
$\frac{n(n-1)}2$ identities
$$
\beta_{(r,s)}\sigma_{(r,s)} = {\alpha_{(r,s)}}, \ 0 \leq r \leq n,
2 \leq s \leq n, r+s \leq n, \tag3.2.10
$$
are defining relations for the ring
 $\Bbb Q[z_1,\overline z_1,\dots, z_n,\overline z_n]^{S_n}$ of
bisymmetric functions (over the rationals).
\endproclaim

\demo{Proof} Since $\Bbb Q[z_1,\overline z_1,\dots, z_n,\overline
z_n]^{S_n}$ is a domain, it injects into its quotient field. In
the quotient field, the identities (3.2.10) are equivalent to the
identities (3.2.8). \qed
\enddemo

At the risk of making a mountain out of a molehill we note that
the complexification $\Bbb R [T^{\Bbb C}\big/S_n]_{\Bbb C}$ of the
real coordinate ring $\Bbb R [T^{\Bbb C}\big/S_n]$ which we are
looking for is simply the smallest subring of the field $\Bbb
Q(z_1,\overline z_1,\dots, z_n,\overline z_n)^{S_n}$ of
bisymmetric functions containing $\Bbb Q[z_1,\overline z_1,\dots,
z_n,\overline z_n]^{S_n}$  and the inverses $\sigma_n^{-1}$ and
$\overline \sigma_n^{-1}$.

\smallskip
\noindent {\smc (3.3) The real semialgebraic structure of the
adjoint quotient.\/}

\noindent
On $\Bbb C^n \times {\Bbb C}^n$, consider the complex conjugation
given by $(\bold z_1, {\bold z_2})
\mapsto (\overline{\bold z}_2, \overline{\bold z}_1)$ ($(\bold
z_1,\bold z_2) \in \Bbb C^n \times \Bbb C^n$).
The real linear subspace invariant under this conjugation is the diagonal space
$$
\Delta=\{(\bold z, \overline{\bold z}); \bold z \in \Bbb C^n\}
\subseteq \Bbb C^n \times {\Bbb C}^n.
$$
Relative to the obvious diagoinal $S_n$-action on
$\Bbb C^n \times {\Bbb C}^n$, $\Delta$ is an
$S_n$-invariant subspace. In terms of the
real coordinates $x_j$ and $y_j$ given by $z_j=x_j+ i y_j$ ($1
\leq j \leq n$), $\Delta$ is the real vector space with
coordinates $x_1$, $y_1$, \dots, $x_n$, $y_n$ and, as a real
non-singular variety, the complex torus $T^{\Bbb C}$ lies in
$\Delta$ as the subspace of points $(x_1,y_1,\dots,x_n,y_n)$ with
$x_j+ iy_j \ne 0$ for $1 \leq j \leq n$. The observations spelled
out in Section 2 reduce the description of the structure of the
adjoint quotient $T^{\Bbb C}/S_n$ to
that of the orbit space $\Delta/S_n$. To determine the latter we
note first that, in view of Corollary 3.2.9, the complexification
$\Bbb R[\Delta/S_n]_{\Bbb C}$ of the real coordinate ring $\Bbb
R[\Delta/S_n]$ of the orbit space $\Delta/S_n$ is generated by the
$d=\frac{n(n+3)}2$ invariants (3.2.1), subject to the relations
(3.2.10).

To extract a description of the orbit space $\Delta/S_n$ as a real
semialgebraic space, take the invariants (3.2.1) as coordinates on
$\Bbb C^d$ and denote by $Y \subseteq \Bbb C^d$ the complex affine
variety given by the equations (3.2.10). For $j\neq k$ and $j+k
\leq n$, $\sigma_{j,k}$ and $\sigma_{k,j}$ are here considered as
independent complex variables. Let
$$
p=(\sigma_{(1,0)}, \sigma_{(0,1)}, \sigma_{(2,0)},
\sigma_{(1,1)},\dots, \sigma_{(0,n)})\colon \Bbb C^n \times \Bbb
C^n \longrightarrow \Bbb C^d \tag3.3.1
$$
be the resulting {\smc Hilbert\/} map, so that the complex variety
$Y \subseteq \Bbb C^d$ coincides with the image of $p$.

The complex conjugation  on $\Bbb C^n \times  {\Bbb C}^n$
introduced above induces the complex conjugation on $\Bbb C^d$
which sends $\sigma_{(j,k)}$ to $\sigma_{(k,j)}$ ($1 \leq j+k \leq n$)
and $\sigma_{(\ell,\ell)}$ ($2 \leq 2\ell \leq n$) to its complex
conjugate $\overline \sigma_{(\ell,\ell)}$. Relative to this complex
conjugation, embed $\Bbb R^d$ into $\Bbb C^d$ in the obvious
manner, that is, as the subspace of invariants under this complex
conjugation. Then the {\it real\/} categorical quotient $\widehat
{\Delta/S_n}$ amounts to the intersection
$Y \cap \Bbb R^d$,
and the actual orbit space ${\Delta/S_n}$ lies in $\widehat
{\Delta/S_n}$ as the semialgebraic set $p(\Delta)$ in $Y \cap
\Bbb R^d$. We will now explain how inequalities defining this
semialgebraic set may be derived; such inequalities are more easily
spelled out in terms of the invariants (3.2.2).

For convenience we will consider the elements of the target
$ \Bbb C^{d}$ as column vectors.
The invariants (3.3.2) furnish the Hilbert map
$$
\Phi=\left[\tau_{(1,0)},
\tau_{(0,1)},\tau_{(2,0)},\tau_{(1,1)},\tau_{(0,2)},\ldots,
\tau_{(0,n)}\right]^{\roman t} \colon \Bbb C^n \times \Bbb C^n
\longrightarrow \Bbb C^{d} \tag3.3.2
$$
where $[\ldots]^{\roman t}$ refers to the transpose of the row
vector $[\ldots]$. In the chosen coordinates, the Jacobian
$J_{\tau}$ of $\Phi$ is the $(d \times (2n))$-matrix
$$
J_{\Phi} = \left[\frac {\partial \Phi}{\partial z_1}, \frac
{\partial \Phi}{\partial {\overline z}_1}, \ldots, \frac {\partial
\Phi}{\partial z_n}, \frac {\partial \Phi}{\partial {\overline
z}_n}\right], \tag3.3.3
$$
with column vectors $\frac{\partial \Phi}{\partial z_j}$ and
$\frac {\partial \Phi}{\partial {\overline z}_j}$, for $1\leq j
\leq n$. A result in \cite\procschw, see also Theorem 6.2 in
\cite\gwschwat, implies at once the following where $B^{\roman t}$
refers to the transpose of the matrix $B$.

\proclaim{Proposition 3.3.4} Inequalities defining the quotient
$\Delta/S_n$ as a real semi-algebraic space in $\Bbb R^{d}$ are
given by the requirement that the hermitian $(d\times d)$-matrix $
J_{\Phi} \overline {J}_{\Phi}^{\roman t} $ be positive
semidefinite. \qed
\endproclaim

Notice that the condition that a given hermitian matrix $B$ be
positive semidefinite is equivalent to a simultaneous system of
inequalities $\{B_{\alpha} \geq 0\}$ where $B_{\alpha}$ runs over
the determinants of the hermitian minors of $B$.

\smallskip
\noindent {\smc (3.4) The stratified symplectic Poisson structure
of the adjoint quotient.\/}
For convenience, on $\roman {GL}(n,\Bbb C)$,
we will take the K\"ahler potential coming from
the invariant quadratic form
$$
 \frak t \to \Bbb R,\
(it_1,\dots, i t_n) \longmapsto t_1^2 + \dots + t_n^2
\tag3.4.1
$$
on the Lie algebra $\frak t$ of the standard maximal torus
of $\roman U(n)$.

\proclaim{Theorem 3.4.2} The Poisson brackets among the
multiplicative generators $\tau_{(j,k)}$ {\rm ($0\leq j \leq n, 0
\leq k \leq n, 1 \leq j+k \leq n$)} of the ring $\Bbb
Q[z_1,\dots,z_n,\overline z_1,\dots,\overline z_n]^{S_n}$ of
bisymmetric functions are given by the formulas
$$
\frac i2\{\tau_{(j_1,k_1)},\tau_{(j_2,k_2)}\} =
(j_1k_2-j_2k_1)\tau_{(j_1+j_2,k_1+k_2)}. \tag3.4.3
$$
\endproclaim

In the statement of this theorem, when $j_1+j_2+k_1+k_2 >n$, the
right-hand side of (3.4.3) is to be rewritten as a polynomial in
terms of the multiplicative generators $\tau_{(j,k)}$, for $0\leq
j \leq n, 0 \leq k \leq n, 1 \leq j+k \leq n$. The theorem
justifies the claim made earlier that the stratified symplectic
Poisson structure is most conveniently described in terms of the
power sum bisymmetric functions.

\demo{Proof} The assertion is established by a straightforward
calculation relying on the formula (2.1.6). We leave the details
to the reader. \qed
\enddemo

Theorem 3.4.2 yields the stratified symplectic Poisson algebra of the
adjoint quotient for the unitary group $\roman U(n)$;
for the group $\roman {SU}(n)$ a slight modification is necessary which
we will explain shortly. The theorem has
the following attractive consequence, the proof
of which is immediate.

\proclaim{Corollary 3.4.4} As a Poisson algebra, the ring $\Bbb
Q[z_1,\dots,z_n,\overline z_1,\dots,\overline z_n]^{S_n}$ of
bisymmetric functions is generated by the power sum functions
$\tau_1$, \dots, $\tau_n$, $\overline \tau_1$, \dots, $\overline
\tau_n$. That is to say:
The ring of bisymmetric functions is
generated by the power sum functions $\tau_1$, \dots, $\tau_n$,
$\overline \tau_1$, \dots, $\overline \tau_n$ together with
iterated Poisson brackets in these functions.
Consequently the ring of bisymmetric functions is
generated as well by the  elementary symmetric
functions $\sigma_1$, \dots, $\sigma_n$,
$\overline \sigma_1$, \dots, $\overline \sigma_n$ together with
iterated Poisson brackets in these functions. \qed
\endproclaim

This corollary establishes the theorem in the introduction
for the special case where $K= \roman{U}(n)$.

\noindent {\smc Remark 3.4.5.\/} There is an intimate relationship
between the rank of the Poisson structure on the adjoint quotient
and the real semi-algebraic structure. Indeed, let $P$ be the
$((2n)\times(2n))$-matrix
$$
P= \left[ \matrix
0 & \{z_1,\overline z_1\} &   0 &   0 & \ldots & 0 & 0
\\
\{\overline z_1, z_1\} & 0&   0 &   0 & \ldots & 0 & 0
\\
\ldot &              \ldot&\ldot&\ldot& \ldots & 0 & 0
\\
\ldot &              \ldot&\ldot&\ldot& \ldots & 0
&\{z_n,\overline z_n\}
\\
\ldot &              \ldot&\ldot&\ldot& \ldots & \{\overline
z_n,z_n\} & 0
\endmatrix
\right] .
$$
The skew-symmetric $(d\times d)$-matrix
$\left[\{\tau_{(j_1,k_1)},\tau_{(j_2,k_2)}\}\right]$ of Poisson
brackets among the multiplicative generators (3.2.2) may be written as
$$
\left[\{\tau_{(j_1,k_1)},\tau_{(j_2,k_2)}\}\right] = J_{\Phi}
P{J}_{\Phi}^{\roman t}. \tag3.4.6
$$
The adjoint quotient $T^{\Bbb C}\big/S_n$ being realized in $\Bbb
C^d$, given a point $(\bold z,\overline{\bold z})$ of $T^{\Bbb
C}$, at the point $\Phi(\bold z,\overline{\bold z})$ of $T^{\Bbb
C}\big/S_n\subseteq \Bbb C^d$, the rank of the Poisson structure
is given by the rank of the matrix (3.4.6) at $\Phi(\bold
z,\overline{\bold z})$; this rank coincides with the rank of the
matrix $J_{\Phi}(\bold z,\overline{\bold z})$ since, at the point
$(\bold z,\overline{\bold z})$ of $T^{\Bbb C}$, the matrix $P$ has
maximal rank. Consequently the rank of the Poisson structure at
the point $\Phi(\bold z,\overline{\bold z})$ of $T^{\Bbb
C}\big/S_n\subseteq \Bbb C^d$ coincides with the rank of the $(d
\times d)$- matrix $J_{\Phi} \overline {J}_{\Phi}^{\roman t}$
coming into play in (3.3) above. This implies that certain {\it
unions of closures of  strata may be described by requiring that
some of the inequalities giving the real semi-algebraic structure
be equalities.\/} Given a stratum $Y$ of $T^{\Bbb C}\big/S_n$,
when this stratum has dimension $2k$, as a smooth symplectic
manifold the stratum arises as the base of a regular finite
covering having as total space a smooth $(2k)$-dimensional
symplectic submanifold of $T^{\Bbb C}$ whence the Poisson bracket,
restricted to the stratum $Y$, necessarily has rank $2k$.

\smallskip
\noindent (3.5) The K\"ahler potential (1.2), restricted to
$T^{\Bbb C}$, can be written as a real analytic function of the
invariants (3.2.1) or (3.2.2) and hence descends to a real
analytic function $\kappa_{\roman {red}}$ on the adjoint quotient,
that is, determines an element of $C^{\omega}(T^{\Bbb
C}\big/S_n)$. We shall illustrate this for $K=\roman {SU}(2)$
below.
\smallskip
\noindent (3.6) $K=\roman {U}(2)$; now the map (3.1.1) (for $n=2$)
induces a complex algebraic isomorphism from the adjoint quotient
onto the complex manifold $\Bbb C \times \Bbb C^*$, with complex
coordinates $\sigma_1$ and $\sigma_2$, $\sigma_2 \ne 0$. There is
one lower stratum which corresponds to the center of $\roman U(2)$
and arises from the diagonal in $T^{\Bbb C} \cong \Bbb C^* \times
\Bbb C^*$. As an affine complex variety, this stratum is the {\it
discriminant variety\/}, and an equation for it arises from
rewriting the $S_2$-invariant equation $(z_1-z_2)^2 = 0$ in terms
of the $S_2$-invariants $\sigma_1$ and $\sigma_2$; this yields the
familiar equation
$$
\sigma_1^2 -4 \sigma_2 = 0\ (\sigma_2 \ne 0),\tag3.6.1
$$
the left-hand side of this equation being the discriminant of the
polynomial $P(z) = z^2 -\sigma_1 z + \sigma_2$. Thus geometrically
the lower stratum amounts to a (complex) parabola in the adjoint
quotient $\Bbb C \times \Bbb C^*$, that is, to a parabola with the
origin removed.

To describe the Poisson structure on the complexification $\Bbb
R[T^{\Bbb C}/S_2]_{\Bbb C}$ of the real coordinate ring $\Bbb
R[T^{\Bbb C}/S_2]$ of the adjoint quotient $T^{\Bbb C}/S_2$,
we note first that, in
view of the observations spelled out in Subsection 3.2 above,
 the algebra of $S_2$-invariants $\Bbb
C[z_1,\overline z_1,z_2,\overline z_2]^{S_2}$ is generated by the
elementary bisymmetric functions
$$
\sigma_1 = z_1+z_2,\ \overline \sigma_1 = \overline z_1+\overline
z_2,\ \sigma_2 = z_1z_2,\ \overline \sigma_2 = \overline
z_1\overline z_2,\ \sigma =\sigma_{(1,1)} =z_1  \overline z_2 +
z_2\overline z_1 ;\tag3.6.2
$$
moreover, these generators are subject to the single defining
relation
$$
(\sigma_1^2-4 \sigma_2)(\overline \sigma_1^2-4 \overline \sigma_2)
= (\sigma_1 \overline \sigma_1 -2\sigma)^2. \tag3.6.3
$$
Indeed, in view of the formula (3.2.7),
$$
\overline \sigma_2 = \overline z_1 \overline z_2 =
-\frac{(z_1\overline \sigma_1 - \sigma)(z_2\overline \sigma_1 -
\sigma)} {(z_1-z_2)^2}. \tag3.6.4
$$
A  calculation shows that this identity is equivalent to (3.6.3).
In terms of the five bisymmetric power sum functions $\tau_1 =
\tau_{1,0}$, $\overline \tau_1 = \tau_{0,1}$, $\tau = \tau_{1,1}$,
$\tau_2 = \tau_{2,0}$, $\overline \tau_2 = \tau_{0,2}$, the
(algebraic) Poisson structure is given by Theorem 3.4.2, and a
straightforward calculation yields expressions for the Poisson
structure in terms of the generators (3.6.2).

Under the present circumstances, other calculations yield
$$
J_{\Phi} \overline {J}_{\Phi}^{\roman t} = \left[\matrix
2 & \overline \tau_1 &0 & 0 &\tau_1
\\
\tau_1 & \tau &0 & 0 &\tau_2
\\
0 & 0 & 2 & \tau_1 & \overline \tau_1
\\
0 & 0 & \overline \tau_1 & \tau & \overline \tau_2
\\
\overline \tau_1 & \overline \tau_2 &\tau_1 &\tau_2 & 2 \tau
\endmatrix
\right] \tag3.6.5
$$
and
$$
\left|J_{\Phi} \overline {J}_{\Phi}^{\roman t}\right| = 16(2\tau -
\tau_1 \overline \tau_1) ((2\tau - \tau_1 \overline \tau_1)^2 -
(2\tau_2-\tau_1^2) (2\overline \tau_2-\overline \tau_1^2)) .
\tag3.6.6
$$
Furthermore, straightforward calculation yields
$$
\align
2\tau - \tau_1 \overline \tau_1&= \sigma_1\overline
\sigma_1 - 2 \sigma =(z_1-z_2)(\overline z_1 - \overline z_2)
\\
2\tau_2-\tau_1^2 &= \sigma_1^2 -4 \sigma_2
\\
2\overline \tau_2-\overline \tau_1^2 &= \overline \sigma_1^2 -4
\overline \sigma_2
\\
(2\tau - \tau_1 \overline \tau_1)^2 - (2\tau_2-\tau_1^2)
(2\overline \tau_2-\overline \tau_1^2) &= (\sigma_1\overline
\sigma_1 - 2 \sigma)^2 -(\sigma_1^2 -4 \sigma_2)(\overline
\sigma_1^2 -4 \overline \sigma_2)
\endalign
$$
whence
$$
\left|J_{\Phi} \overline {J}_{\Phi}^{\roman t}\right| = 16
(\sigma_1\overline \sigma_1 - 2 \sigma) ((\sigma_1\overline
\sigma_1 - 2 \sigma)^2 -(\sigma_1^2 -4 \sigma_2)(\overline
\sigma_1^2 -4 \overline \sigma_2)). \tag3.6.7
$$
In view of the
relation (3.6.3), the second factor
on the right-hand side of this equation
is identically zero. However,
the lower stratum is characterized by the equation
$$
\sigma_1\overline \sigma_1 - 2 \sigma = 0 \tag3.6.8
$$
whence we
see that, on the lower stratum, the Poisson structure necessarily
has rank 2. By Proposition 3.3.4, inequalities characterizing the
real semialgebraic structure are given by requiring that, for
$\ell = 1,2,3$, the symmetric $(\ell \times \ell)$-minors of
(3.6.5) be non-negative. This leads to (not necessarily
independent) inequalities of the kind
$$
\alignat1
\tau &\geq 0
\tag3.6.9.1
\\
\left|\matrix 2 & \overline \tau_1
\\
\tau_1 & \tau
\endmatrix \right|
= 2\tau - \tau_1 \overline \tau_1= \sigma_1\overline \sigma_1 - 2
\sigma &\geq 0 \tag3.6.9.2
\\
\left|\matrix 2 & \tau_1& \overline \tau_1
\\
\overline \tau_1 & \tau & \overline \tau_2
\\
\tau_1 & \tau_2 & 2 \tau
\endmatrix \right|
= 4\tau - 2\tau_2 \overline \tau_2 -3 \tau \tau_1 \overline \tau_1
+\tau_1^2\overline \tau_2 +  \overline\tau_1^2 \tau_2&\geq 0
\tag3.6.9.3
\\
\left|\matrix \tau & 0& 0 & \tau_2
\\
0& 2 & \tau_1 & \overline \tau_1
\\
0&\overline \tau_1 & \tau &  \overline \tau_2
\\
\overline \tau_2 &\tau_1 & \tau_2 & 2 \tau
\endmatrix \right|
= (\tau^2 - \tau_2 \overline \tau_2)(2 \tau - \tau_1 \overline
\tau_1) &\geq 0 \tag3.6.9.4
\endalignat
$$
etc. We note that, in view of the defining relation (3.6.3),
$$
\align \tau^2 - \tau_2 \overline \tau_2 &= \frac 14
((\sigma_1\overline \sigma_1 - 2 \sigma)^2 -(\sigma_1^2 -4
\sigma_2)(\overline \sigma_1^2 -4 \overline \sigma_2)) + \sigma_2
\overline \sigma_1^2 + \overline \sigma_2 \sigma_1^2 - \sigma_1
\overline \sigma_1 \sigma
\\
&= \sigma_2 \overline \sigma_1^2 + \overline \sigma_2 \sigma_1^2 -
\sigma_1 \overline \sigma_1 \sigma = 4 \sigma_2 \overline \sigma_2
- \sigma^2
\endalign
$$
whence the inequality (3.6.9.4) is equivalent to
$$
(4 \sigma_2 \overline \sigma_2 - \sigma^2) (\sigma_1\overline
\sigma_1 - 2 \sigma) \geq 0. \tag3.6.9.5
$$

To elucidate the real structure of the quotient, in particular its
semialgebraicity, entirely in terms of appropriate real data, let
$z_1=x_1+iy_1$, $z_2=x_2+iy_2$, and write $\sigma_1 = X+iY$ and
$\sigma_2 = U+iV$, so that
$$
\gathered X= x_1 + x_2, \ Y = y_1+y_2,\ U= x_1x_2 -y_1y_2,\ V=x_1
y_2 +x_2 y_1,
\\
\sigma =2(x_1 x_2 +y_1 y_2).
\endgathered
 \tag3.6.10
$$
In view of the observations spelled out in Section 3 above,
the algebra of $S_2$-invariants $\Bbb Q[x_1, y_1,x_2, y_2]^{S_2}$
is generated by $X$, $Y$, $U$, $V$, and $\sigma$, and a
straightforward calculation shows that the defining relation
(3.6.3), written out in terms of these generators, has the form
$$
(X^2 - Y^2 -4U)^2 + 4(XY-2V)^2 = (X^2 + Y^2 -2 \sigma)^2.
\tag3.6.11
$$
Thus the {\it real categorical\/} quotient of $T^{\Bbb C} \cong
\Bbb C^* \times \Bbb C^*$ modulo $W=S_2$ is realized in $\Bbb R^5$
with coordinates $X,Y,U,V,\sigma$ as the real variety given by the
equation (3.6.11) with $U \ne 0$ and, since $\sigma_1 \overline
\sigma_1 = X^2 + Y^2$, in view of the inequalities (3.6.9.2) and
(3.6.9.5), the quotient $T^{\Bbb C}\big /W$ we are really
interested in is the semialgebraic subspace of the real
categorical quotient given by the inequalities
$$
\align X^2 + Y^2 - 2\sigma &\geq 0
\tag3.6.14
\\
4(U^2 + V^2) - \sigma^2 &\geq 0. \tag3.6.15
\endalign
$$
It is straightforward to rewrite the Poisson brackets (3.4.3) in
terms of the generators (3.6.2) and the generators (3.6.10). For
example
$$
\align
\frac i2
\{\sigma_1,\overline \sigma_1\} &=
\{X,Y\} = \tau = X^2+Y^2-\sigma,
\\
\frac i2
\{\sigma_2,\overline \sigma_2\} &=
\{U,V\} = 2\sigma_2 \overline \sigma_2 = 2(U^2+V^2),
\\
\frac i2
\{\sigma_1,\overline \sigma_2\} &=
\sigma_1 \overline \sigma_2 ,
\endalign
$$
etc. We leave the details to the reader.
\smallskip
\noindent (3.7) $K=\roman {U}(3)$; in this case the map (3.1.1)
(for $n=3$) induces a complex algebraic isomorphism from the
adjoint quotient $T^{\Bbb C}\big/S_3 \cong (\Bbb C^* \times \Bbb
C^* \times \Bbb C^*)\big/S_3$ onto the complex
algebraic
manifold $\Bbb C
\times \Bbb C\times \Bbb C^*$, with complex coordinates
$\sigma_1$, $\sigma_2$, and $\sigma_3$, $\sigma_3 \ne 0$. The
adjoint quotient has three strata; in terms of the maximal torus
$T^{\Bbb C} \cong \Bbb C^* \times \Bbb C^*\times \Bbb C^*$, these
arise from matrices $\roman{diag}(z_1,z_2,z_3)\in T^{\Bbb C}$ of
the kind $z_1\ne z_2\ne z_3\ne z_1$, $z_1= z_2\ne z_3$, and $z_1=
z_2=z_3$. The complement of the 3-dimensional stratum is the
corresponding {\it discriminant variety\/}, that is, the affine
variety given by the equation
$$
D_3(1,-\sigma_1,\sigma_2,-\sigma_3) = 0
$$
where
$$
D_3(a_0,a_1,a_2,a_3)
=a_1^2a_2^2 -4 a_0a_2^3 -4 a_1^3 a_3 -27 a_0^2a_3^2 + 18 a_0 a_1
a_2 a_3
\tag3.7.1
$$
refers to the discriminant of the polynomial
$a_0 w^3 + a_1 w^2 + a_2 w + a_3$.
The discriminant variety is the closure of the 2-dimensional stratum.
Likewise, when
$$
D_2(b_0,b_1,b_2) = b_1^2 - 4 b_0 b_2
\tag3.7.2
$$
refers to the discriminant of the polynomial
$b_0 w^2 + b_1 w + b_2$,
the 1-dimensional stratum
is the affine curve given by the additional equation
$$
D_2(3,-2\sigma_1,\sigma_2) = 0.
$$
This is just the cubic complex curve given by the parametrization
$v \mapsto (3v, 3 v^2, v^3)$ ($v \in \Bbb C$).

In view of the observations spelled out in (3.2) above, the
complexification $\Bbb R[T^{\Bbb C}/S_3]_{\Bbb C}$ of the real
coordinate ring $\Bbb R[T^{\Bbb C}/S_3]$ of the adjoint quotient
$T^{\Bbb C}/S_3$ is the algebra of $S_3$-invariants $\Bbb C[z_1,\overline
z_1,z_2,\overline z_2,z_3,\overline z_3]^{S_3}$, and this algebra
is generated by the elementary bisymmetric functions
$$
\gathered \sigma_1= \sigma_{(1,0)},
\ \overline \sigma_1= \sigma_{(0,1)},
\ \sigma_2= \sigma_{(2,0)},\ \overline
\sigma_2= \sigma_{(0,2)}, \ \sigma_3= \sigma_{(3,0)},
\ \overline \sigma_3= \sigma_{(0,3)} ,
\\
\sigma = \sigma_{(1,1)}
=z_1 \overline z_2 + z_2\overline z_1 + z_2 \overline z_3 +
z_3\overline z_2 + z_3 \overline z_1 + z_1\overline z_3,
\\
\rho= \sigma_{(2,1)}=
(z_1\overline z_2 z_3 + z_2 \overline z_3 z_1 + z_3\overline
z_1 z_2),
\\
\overline \rho = \sigma_{(1,2)}= (z_1\overline z_2 \overline z_3 +
z_2 \overline z_3 \overline z_1 + z_3\overline z_1 \overline z_2).
\endgathered
\tag3.7.3
$$
Moreover, in view of Proposition 3.2.3,
 the field of bisymmetric functions is freely generated by
$\sigma_1$, $\overline \sigma_1$,
$\sigma_2$, $\sigma_3$, $\sigma$ and $\rho$.

In the present special case where $n=3$, the formula (3.2.7) comes
down to
$$
\overline z_1 = \frac
{z^2_1\overline \sigma_1 - z_1\sigma + \rho}
{(z_1-z_2)(z_1-z_3)},
\
\overline z_2 = \frac
{z^2_2\overline \sigma_1 - z_2\sigma + \rho}
{(z_2-z_1)(z_2-z_3)},
\
\overline z_3 = \frac
{z^2_3\overline \sigma_1 - z_3\sigma + \rho}
{(z_3-z_1)(z_3-z_2)} .
$$
Thus the symmetric function
$
\overline \sigma_2 = \overline z_1 \overline z_2 +
\overline z_1 \overline z_3 +
\overline z_2\overline z_3
$
takes the form
$$
\overline \sigma_2 = \frac{\sum_{j \roman{mod}
3}(z_{j+2}-z_j)(z_{j+2}-z_{j+1}) {(z_j^2\overline \sigma_1 -
z_j\sigma + \rho)} {(z_{j+1}^2\overline \sigma_1 - z_{j+1}\sigma +
\rho)}} {D_3(1,-\sigma_1,\sigma_2, - \sigma_3)}, \tag3.7.4
$$
and this identity yields a relation of the kind
$$
D_3(1,-\sigma_1,\sigma_2, - \sigma_3) \overline \sigma_2
=\alpha(\sigma_1,\sigma_2,\sigma_3,\overline\sigma_1, \sigma,
\rho) \tag3.7.5
$$
where
$\alpha(\sigma_1,\sigma_2,\sigma_3,\overline\sigma_1, \sigma, \rho)$
is a polynomial in the indicated generators.
A tedious but straightforward
calculation involving the Newton polynomials
yields the following expression for $\alpha$:
$$
\aligned
\alpha(\sigma_1,\sigma_2,\sigma_3,\overline\sigma_1, \sigma, \rho)
&=(9 \sigma_3^2 + \sigma_2^3 - 4 \sigma_1 \sigma_2\sigma_3)
\overline \sigma_1^2
\\
&
\quad
+(4 \sigma_1^2 \sigma_3 -3 \sigma_2 \sigma_3 - \sigma_1 \sigma_2^2)
\overline \sigma_1 \sigma
\\
&
\quad
+(6 \sigma_1 \sigma_3 - \sigma_2^2) \overline \sigma_1 \rho
+(\sigma_2^2-3\sigma_1 \sigma_3 )  \sigma^2
\\
&
\quad
+(9 \sigma_3 - \sigma_1 \sigma_2) \sigma \rho
+(\sigma_1^2 - 3 \sigma_2) \rho^2
\endaligned
\tag3.7.6
$$
Likewise, the symmetric function
$
\overline \sigma_3 = \overline z_1 \overline z_2 \overline z_3
$
takes the form
$$
\overline \sigma_3 = \frac{(z^2_1\overline \sigma_1 - z_1\sigma +
\rho) (z^2_2\overline \sigma_1 - z_2\sigma + \rho) (z^2_3\overline
\sigma_1 - z_3\sigma + \rho)} {D_3(1,-\sigma_1,\sigma_2, -
\sigma_3)}, \tag3.7.7
$$
and this identity yields a relation of the kind
$$
D_3(1,-\sigma_1,\sigma_2, - \sigma_3) \overline \sigma_3
=\beta(\sigma_1,\sigma_2,\sigma_3,\overline\sigma_1, \sigma, \rho)
\tag3.7.8
$$
where
$\beta(\sigma_1,\sigma_2,\sigma_3,\overline\sigma_1, \sigma, \rho)$
is a polynomial in the indicated generators.
Again
a tedious but straightforward
calculation yields the following expression for $\beta$:
$$
\aligned
\beta(\sigma_1,\sigma_2,\sigma_3,\overline\sigma_1, \sigma, \rho)
&=
\sigma_3^2 \overline \sigma_1^3
-\sigma_2 \sigma_3 \overline \sigma_1^2 \sigma
+
(\sigma_2^2 - 2 \sigma_1 \sigma_3) \overline \sigma_1^2 \rho
+ \sigma_1 \sigma_3 \overline \sigma_1 \sigma^2
\\
&\quad
-((\sigma_1^2 - 2 \sigma_2)\sigma_1
- \sigma_1^3 + 3 \sigma_1 \sigma_2 - 3 \sigma_3)
\overline \sigma_1 \sigma \rho
\\
&\quad
-\sigma_3 \sigma^3
+(\sigma_1^2 - 2 \sigma_2) \overline \sigma_1 \rho^2
\\
&\quad+ \sigma_2 \sigma^2 \rho
- \sigma_1 \sigma \rho^2
+ \rho^3
\endaligned
\tag3.7.9
$$
Moreover, a calculation yields the relation
$$
(\sigma_1^2-4 \sigma_2)(\overline \sigma_1^2-4 \overline \sigma_2)
= (\sigma_1 \overline \sigma_1 -2\sigma)^2 + 2\rho \overline
\sigma_1 + 2\overline\rho \sigma_1. \tag3.7.10
$$
Since substituting the right-hand side of (3.7.4) for
$\overline \sigma_2$ in (3.7.10) yields an identity which, in turn,
yields $\overline \rho = \sigma_{(1,2)}$ as a rational function of
the generators $\sigma_1$, $\sigma_2$, $\sigma_3$,
$\overline\sigma_1$, $\sigma$, $\rho$, we conclude that, by virtue
of Corollary 3.2.9, the relations (3.7.5), (3.7.8), and (3.7.10)
are {\it defining\/} relations for the algebra $\Bbb
C[z_1,\overline z_1,z_2,\overline z_2,z_3,\overline z_3]^{S_3}$ of
$S_3$-invariants. We note that when we set the variables $z_3$ and
$\overline z_3$ equal to zero, that is, under the obvious
surjection
$$
\Bbb C[z_1,\overline z_1,z_2,\overline z_2,z_3,\overline
z_3]^{S_3} @>>> \Bbb C[z_1,\overline z_1,z_2,\overline z_2]^{S_2},
$$
the relation (3.7.10) passes to the relation (3.6.3).

Straightforward calculation yields the Poisson structure given by
Theorem 3.4.2 in terms of the generators (3.7.3) and, furthermore,
the corresponding real semialgebraic structure made explicit in
(3.3) above in terms of appropriate equations and inequalities,
cf. Proposition 3.3.4. We refrain from spelling out details.

\smallskip
\noindent (3.8) $K=\roman {SU}(n)$. The complexification $\roman
{SU}(n)^{\Bbb C}$ equals the group $\roman {SL}(n,\Bbb C)$, the
standard maximal torus $T^{\Bbb C} \cong (\Bbb C^*)^{n-1}$
consists of the diagonal matrices in $\roman {SL}(n,\Bbb C)$, that
is, diagonal matrices of the kind $\roman{diag}(w_1,\dots,w_n)$
with $w_1\cdot\dots \cdot w_n =1$,
 and the Weyl group $W$ is the
symmetric group $S_n$ on $n$ letters which acts on $T^{\Bbb C}$ by
permutation of the diagonal entries $w_1,\dots,w_n$ of the
matrices $\roman{diag}(w_1,\dots,w_n)$. Hence
 the complex affine coordinate ring $\Bbb C[T^{\Bbb
C}]$ may be written as
$$
\Bbb C[T^{\Bbb C}] \cong \Bbb C[v_1,\dots, v_n]\big/ (\sigma_n -1)
$$
in such a way that the $W$-action permutes the $v_1$, \dots,
$v_n$. Thus we can carry out the construction in (2.2) above with
$m=n$ and $V= \Bbb C^n$ the $S_n$-representation where
$S_n$ acts on $\Bbb C^n$ by permutation of the factors $\Bbb
C$ so that complex algebraically the quotient $V/S_n$ amounts to
the $n$'th symmetric power of a copy of $\Bbb C$. This
construction yields $T^{\Bbb C}$ as the complex $S_n$-subvariety of
$\Bbb C^n$ given by the single equation
$$
v_1\cdot \dots \cdot v_n= 1.
$$
This observation, in turn, implies at once that the adjoint quotient
$T^{\Bbb C}\big/S_n$
 may be realized by means of the map
$$
(\sigma_1,\dots,\sigma_{n-1}) \colon T^{\Bbb C} @>>> \Bbb C^{n-1}.
\tag3.8.1
$$
This description of the categorical quotient of $\roman
{SL}(n,\Bbb C)$ as an affine $(n-1)$-dimensional space is
consistent with Steinberg's observation quoted earlier; indeed,
the elementary symmetric functions $\sigma_1,\dots,\sigma_{n-1}$
are the restrictions to $T^{\Bbb C}$ of the characters of the
fundamental (finite dimensional) representations of
$\roman{SL}(n,\Bbb C)$.

In view of the observations spelled out in Subsection 3.2 above, the
complexification $\Bbb R[T^{\Bbb C}/S_n]_{\Bbb C}$ of the real
coordinate ring $\Bbb R[T^{\Bbb C}/S_n]$ of the adjoint quotient
under discussion is generated by the $\frac{n(n+3)}2$ invariants
of the kind (3.2.1); however, in view of the relations $\sigma_n
=1$ and $\overline\sigma_n =1$, we may discard the generators
$\sigma_n$ and $\overline\sigma_n$. Inequalities characterizing the real
semialgebraic structure may be derived from the corresponding
inequalities for the case where $K= \roman {U}(n)$.

We will now explain how a modification of the statement of Theorem 3.4.2
yields the
stratified symplectic Poisson algebra
on the adjoint quotient $T^{\Bbb C}/S_n$. It is obvious that the generators
(3.2.2) yield multiplicative generators for the complexification of
the real coordinate ring of this quotient as well.
The problem is that we cannot
simply restrict the Poisson brackets to
$\roman{SL}(n,\Bbb C)$, viewed as a subgroup of $\roman{GL}(n,\Bbb C)$
since, in the language of {\it constrained systems\/},
$\roman{SL}(n,\Bbb C)$ is a {\it second class\/} constraint in
$\roman{GL}(n,\Bbb C)$, cf. \cite\diracone.
On $\roman{SL}(n,\Bbb C)$
we take the K\"ahler potential given
by restriction of the K\"ahler potential on
$\roman {GL}(n,\Bbb C)$ determined by (3.4.1) above.
Then the embedding of
$\roman{SL}(n,\Bbb C)$ into $\roman{GL}(n,\Bbb C)$
is one of K\"ahler manifolds.
In view of the discussion in Section 2 above,
the complexification
 $\Bbb R[T^{\Bbb C}/S_n]_{\Bbb C}$ of the real
coordinate ring of the adjoint quotient $T^{\Bbb C}/S_n$
amounts to the ring of invariants
$$
\left(\Bbb C[z_1,\dots,z_n,\overline z_1,\dots,\overline z_n]
\big/(\sigma_n-1,\overline \sigma_n -1)\right)
^{S_n}.
$$
We will write the Poisson bracket on the ring
$\left(\Bbb C[z_1,\dots,z_n,\overline z_1,\dots,\overline z_n]\right)
^{S_n}$
given in Theorem 3.4.2 above as $\{\{\,\cdot\, , \,\cdot\,\}\}$.

\proclaim{Theorem 3.8.2}
On the complexification $\Bbb R[T^{\Bbb C}/S_n]_{\Bbb C}$
of the real coordinate ring of the adjoint quotient $T^{\Bbb C}/S_n$,
the Poisson brackets among the
multiplicative generators
$$
\tau_{(j,k)}, \quad 0\leq j \leq n, 0
\leq k \leq n, 1 \leq j+k \leq n,
$$
are given by the formulas
$$
\aligned
\frac i2\{\tau_{(j_1,k_1)},\tau_{(j_2,k_2)}\} &=
\frac i2\{\{\tau_{(j_1,k_1)},\tau_{(j_2,k_2)}\}\}
\\
& \quad
-
\frac i2
\frac{\{\{\tau_{(j_1,k_1)},\overline \sigma_n\}\}
\{\{\sigma_n,\tau_{(j_2,k_2)}\}\}}
{\{\{\sigma_n,\overline \sigma_n\}\}}
\\
& \quad
-
\frac i2
\frac{\{\{\tau_{(j_1,k_1)}, \sigma_n\}\}
\{\{\overline\sigma_n,\tau_{(j_2,k_2)}\}\}}
{\{\{\overline \sigma_n,\sigma_n\}\}}.
\endaligned
\tag3.8.3
$$
\endproclaim

\demo{Proof} We recall that, in general, given a smooth symplectic
manifold and a smooth symplectic submanifold, the Poisson bracket
on the submanifold is induced from the {\it Dirac bracket\/} on
the ambient manifold, cf. \cite\diracone. Hence, by virtue of
Theorem 3.4.2, the statement of the Theorem is a consequence of
the fact that, given two functions $f$ and $h$ on the
complexification $\widetilde T^{\Bbb C}$ ($\subseteq
\roman{GL}(n,\Bbb C)$) of the standard maximal torus $\widetilde
T$ in $\roman U(n)$, relative to the constraints $\sigma_n=1$ and
$\overline \sigma_n=1$, the {\it Dirac bracket\/} $\{f,h\}$ is
given by
$$
\{f,h\} =\{\{f,h\}\}
-\frac {\{\{f,\sigma_n-1\}\} \{\{\overline \sigma_n-1,h\}\}}
{\{\{\overline \sigma_n-1,\sigma_n-1\}\}}
-\frac {\{\{f,\overline \sigma_n-1\}\} \{\{\sigma_n-1,h\}\}}
{\{\{\sigma_n-1,\overline \sigma_n-1\}\}}. \qed
$$
\enddemo

The following is immediate, cf. Corollary 3.4.4.

\proclaim{Corollary 3.8.4} As a Poisson algebra, the
complexification $\Bbb R[T^{\Bbb C}/S_n]_{\Bbb C}$ of the real
coordinate ring of the adjoint quotient $T^{\Bbb C}/S_n$ for
$\roman{SU}(n)$ is generated by the power sum functions $\tau_1$,
\dots, $\tau_n$, $\overline \tau_1$, \dots, $\overline \tau_n$
and hence by the elementary symmetric functions  $\sigma_1$,
\dots, $\sigma_{n-1}$, $\overline \sigma_1$, \dots, $\overline
\sigma_{n-1}$. \qed
\endproclaim

Since the elementary symmetric functions
are the characters of the fundamental representations
of the group $\roman{SU}(n)$,
the statement of the theorem in the introduction
for $K=\roman{SU}(n)$
is an immediate
consequence of this corollary.

\noindent (3.9) $K=\roman {SU}(2)$: We will now make the
stratified K\"ahler structure explicit for $n=2$, and we will
spell out the reduced K\"ahler potential.  The standard maximal
torus $T^{\Bbb C}$ in $K^{\Bbb C} =\roman {SL}(2,\Bbb C)$ is the
group of diagonal matrices of the kind $\left[\matrix z & 0
\\ 0 & z^{-1}
\endmatrix\right]$ ($z \in \Bbb C, z \ne 0$)
in $\roman {SL}(2,\Bbb C)$. Under the identification of $\Bbb C^*$
with $ T^{\Bbb C}$ given by the assignment to $z\in \Bbb C^*$ of
the diagonal matrix $\roman{diag}(z,z^{-1})$, the action of the
Weyl group $W = S_2$ on $T^{\Bbb C}$ amounts to the $S_2$-action
on $\Bbb C^*$ where the non-trivial involution acts by inversion.
The first elementary symmetric function $\sigma_1$, that is, the
character of the defining representation, comes down to the
holomorphic function
$$
\sigma_1=\chi \colon \Bbb C^* @>>> \Bbb C, \quad \chi(z) = z +
z^{-1},\quad z \in \Bbb C,
$$
which identifies the quotient $T^{\Bbb C}\big / S_2 \cong \Bbb
C^*\big / S_2$ with the affine complex line $\Bbb C$.

Next we elucidate the real structure. According to the general
construction spelled out in (3.8), the
complexification $\Bbb R[T^{\Bbb C}/S_n]_{\Bbb C}$ of the real
coordinate ring $\Bbb R[T^{\Bbb C}/S_n]$
is generated by
$Z=\sigma_1 = z + z^{-1}$, $\overline Z =\overline \sigma_1= \overline z +
\overline z^{-1}$, and $\sigma = \frac z{\overline z} + \frac
{\overline z} z$. In view of (3.6.3), these generators are subject
to the relation
$$
\sigma^2 + Z^2 + \overline Z^2 - Z \overline Z \sigma -4 = 0.
\tag3.9.1
$$

In terms of the real and imaginary parts $x$ and
$y$ of the holomorphic coordinate $z = x+iy$, the non-trivial
element of $S_2$ acts on $\Bbb C^*$ via the assignment to $(x,y)$
of $(\frac x {r^2},-\frac y {r^2})$ where as usual $r^2 =x^2 + y^2
= z \overline z$, and the two $S_2$-invariant functions $X =
x+\frac x {r^2}$ and $Y = y-\frac y {r^2}$ on $T^{\Bbb C} \cong
\Bbb C^*$ serve as real coordinate functions on the quotient $\Bbb
C^*\big/S_2$ in such a way that the complex $S_2$-invariant
$Z=\sigma_1= X + i Y =z + z^{-1}$ is a holomorphic coordinate.
Introduce the  $S_2$-invariant function $b$ on $\Bbb C^*$
defined by
$$
2b = r^2 + \frac 1{r^2}. \tag3.9.2
$$
This is a real algebraic
 function on $\Bbb C^*$ which is $S_2$-invariant and hence
descends to the quotient $\Bbb C^*\big/S_2$ as a continuous
function. Since $\sigma = X^2+Y^2 -2b$, the algebra of invariants
$\Bbb R[\Bbb C^*]^{S_2}$ is generated by $X$, $Y$, and $b$, and the relation
$$
Y^2 = b^2 -1 - (b-1)\frac{X^2 + Y^2}2; \tag3.9.3
$$
is equivalent to the relation (3.9.1).
Moreover, on the quotient $\Bbb C^*\big /S_2$, the
function $b$ is non-negative whence this quotient is realized in
$\Bbb R^3$, with coordinates $X$, $Y$, and $b$, as the
semialgebraic set given by the equation (3.9.3) and the obvious
inequality ${b \geq 1}$. Indeed, this inequality may also be
derived from the inequality (3.6.12) and,
since $\sigma = \frac z{\overline z} + \frac
{\overline z} z$, the inequality $\sigma^2
\leq 4$ resulting from (3.6.13) is obviously satisfied.

As for the stratified Poisson structure,
Theorem 3.8.2 or
a direct
calculation yields
$$
\aligned\frac i2\{Z,\overline Z\} &=z \overline z +\frac 1{z \overline z} -
\left(\frac z {\overline z} +\frac {\overline z}z\right)
\\
&=
\tau-\sigma =
Z \overline Z - 2\left(\frac z {\overline z} +\frac {\overline z}
z\right)= Z
\overline Z - 2\sigma.
\endaligned
\tag3.9.4
$$
In terms of the real generators, the formula (3.9.4) comes down to
$$
\{X,Y\} = X^2 + Y^2-2\sigma = 4b -(X^2 + Y^2). \tag3.9.5
$$
Further, it is straightforward to determine the remaining Poisson
brackets among the generators $X$, $Y$, and $b$. Hence, as a
complex algebraic space, the quotient $\Bbb C^*\big/ S_2$ is a
copy of the complex line;  the algebra $\Bbb R[\Bbb C^*\big/ S_2]$
is the algebra of polynomial functions in the variables $X,Y,b$,
subject to the relation (3.9.3); and the stratified symplectic
Poisson bracket on this algebra is given by (3.9.5).

\noindent {\smc 3.9.6. The reduced K\"ahler potential.\/} Recall
that the K\"ahler potential $\kappa$ on $\Bbb C^*$, cf. (1.2) and
(2.1.3), is given by
$$
\kappa(z)=\frac 14 \log^2(z \overline z) = \log^2(r) .
$$
This is a real analytic function which is manifestly
$S_2$-invariant. Hence it descends to a continuous function
$\kappa_{\roman{red}}$  on the
quotient $\Bbb C^*\big/S_2$, and
$\kappa_{\roman{red}}$
 lies in
$C^{\omega}(\Bbb C^*\big/S_2)$;
however, even though
$\Bbb C^*\big/S_2$ is topologically an ordinary plane,
with coordinates $X$ and $Y$,
$\kappa_{\roman{red}}$
is not smooth in the usual sense at the
two singular points and hence is not a smooth function of the two
variables $X$ and $Y$. Indeed, away from the two singular points,
the function $\kappa_{\roman{red}}$ is an ordinary K\"ahler
potential, but at each of these two points the K\"ahler structure
is {\it not even defined\/} as an infinitesimal object on the
{\it ordinary\/} tangent
space to the quotient at each of these points,
relative to the standard smooth structure on
$\Bbb C^*\big/ S_2\cong\Bbb C \cong \Bbb R^2$. Since $Z = z + z^{-1}$, the
coordinate $z$ satisfies the quadratic equation $\left(z-\frac
Z2\right)^2 = \frac{Z^2}4-1$ whence, in view of (1.2), at $Z= X+iY
\ne\pm 2$, the reduced potential $\kappa_{\roman{red}}$ is given
by
$$
\kappa_{\roman{red}}(Z) =\frac 14 \log^2(z \overline z)
=\log^2\left|\frac Z 2 + \sqrt{\frac {Z^2}4 -1}\right| .
\tag3.9.6.1
$$
Notice that, at $Z=\pm 2$, this function is manifestly not an
ordinary smooth function of the variables $X$ and $Y$. To justify
the claim that $\kappa_{\roman{red}}$ lies in $C^{\omega}(\Bbb C^*\big/S_2)$
we must show that $\kappa_{\roman{red}}$ may be written as a real analytic
function of the variables $X$, $Y$ and $b$. In order to see this,
consider the entire analytic function $\roman {Ch}$ in the variable $w$
given by
$$
\roman{Ch}(w) = \sum_{j=0}^{\infty} \frac {w^j}{(2j)!}.
$$
For intelligibility we note that $\roman{Ch}(w^2)=
\roman{cosh}(w)$ and $\roman{Ch}(-w^2)= \roman{cos}(w)$. For $w>
-\pi^2$, the function $\roman{Ch}$ is strictly increasing and
hence admits an analytic inverse
$$
\roman{arCh}\colon ]-1, \infty[ @>>> ]\pi^2, \infty[.
$$
Recall that $t=\log r$ and that $\kappa = t^2$. Since
$$
2b=r^2 + \frac 1{r^2} =2\roman{cosh}(2t) =2\roman{Ch}(4t^2)
=2\roman{Ch}(4\kappa),
$$
for $b \geq 1$, we have
$$
\kappa =\frac 14 \roman{arCh}(b).
$$
A suitable extension of the restriction of the function $
\roman{arCh}$ to $]0, \infty[$ to a smooth function defined on the
entire real line yields $\kappa$ as a smooth function of the
variable $b$ as desired.  Thus $\kappa_{\roman{red}}$ may be
written as a smooth function of the variable $b$, indeed,
as a  {\it real analytic\/} function of the variable $b$ on a neighborhood of
the quotient $\Bbb C^*\big/S_2$, realized in the copy of $\Bbb
R^3$ with coordinates $X$, $Y$, and $b$.

\noindent {\smc 3.9.6.2. Remark.\/} The real Poisson structure
given by (3.9.4) or (3.9.5) {\it cannot\/} be described in terms
of the coordinates $X$ and $Y$ alone. Indeed, the (holomorphic)
derivative $\chi'$ has the two fixed points $1$ and $-1$ of the
$S_2$-action as zeros, and hence $b$ is not an ordinary smooth
function of the variables $X$ and $Y$. Solving (3.9.3) for $b$,
and noticing that $b \geq 1$, we obtain
$$
b = \frac {X^2+Y^2}{4} +\sqrt{1+ \frac {Y^2-X^2} 2+ \frac {(X^2 +
Y^2)^2}{16}}
$$
whence, at $(X,Y)=(\pm 2,0)$, $b$ is not smooth as a function of
the variables $X$ and $Y$. Moreover, ${ b(1) = b(-1) = 1}$ whence
the Poisson bracket $\{X,Y \}$ vanishes  at the two fixed points
of the $S_2$-action.

\noindent {\smc 3.9.7. Geometric Interpretation.\/} The composite
of $\chi$ with the complex exponential mapping
amounts to twice the holomorphic hyperbolic cosine function.
This observation provides a geometric interpretation of the
adjoint quotient under discussion in terms of the geometry
of hyperbolic cosine
and in particular visualizes the familiar fact that, unless there is
a single stratum, the strata arising from cotangent bundle
reduction are {\it not\/} cotangent bundles on strata of the orbit
space of the base space. The real space which underlies
this  quotient has received considerable
attention: It is the \lq\lq canoe\rq\rq, cf. \cite\armcusgo,
\cite\cushbate\ (pp.~148 ff.) and \cite\lermonsj. This space also
arises as the reduced phase space for the spherical pendulum, the
two singular points being absolute equilibria. The resulting complex
analytic interpretation of the spherical pendulum seems to be new, though.
See \cite\bedlepro\ for details.

\smallskip
\noindent (3.10) $K=\roman {SU}(3)$, with maximal torus $T\cong
S^1 \times S^1$. Complex algebraically, the map (3.8.1) for $n=3$
comes down to the map
$$
(\sigma_1,\sigma_2) \colon T^{\Bbb C} \longrightarrow \Bbb C^2,
\tag3.10.1
$$
and this map induces a complex algebraic isomorphism from the
quotient of $\roman {SL}(3,\Bbb C)$, realized as the orbit space
 $T^{\Bbb C}\big/S_3$, onto a copy $\Bbb C^2$ of 2-dimensional
complex affine space, and we take $\sigma_1$ and $\sigma_2$
as complex coordinates on the adjoint quotient. The
discussion in (3.1) and (3.7) above reveals that the complement of
the top stratum amounts to the complex affine curve given by the
equation
$$
D_3(1,-\sigma_1, \sigma_2,-1) = 0,
$$
an explicit expression for
$D_3(1,-\sigma_1,\sigma_2,-1)$
being given by (3.7.1).
This curve is plainly parametrized by the restriction
$$
\Bbb C^* @>>> \Bbb C \times \Bbb C, \quad z \mapsto (2z + z^{-2},
z^2 + 2z^{-1}) \tag3.10.2
$$
of (3.10.1) to the diagonal, and this holomorphic curve
parametrizes the {\it closure\/} of the {\it middle\/} stratum,
that is, of the stratum given by the normalized complex degree 3
polynomials with a proper double root (i.~e. one which is not a
triple root) and with constant coefficient equal to 1. This curve
has the three (complex) singularities $(3,3)$, $3(\eta,\eta^2)$,
$3(\eta^2,\eta)$, where $\eta^3 = 1, \eta \ne 1$. These points are
the images of the (conjugacy classes) of the three central
elements under (3.10.1); as complex curve singularities, these
singularities are cuspidal. These three points constitute the {\it
bottom\/} stratum.

Since, for $z$ with $|z|=1$, $2\overline z + \overline z^{-2}$
equals $z^2+ 2z^{-1}$, the restriction of (3.10.2) to the {\it real\/} torus
$T\subseteq T^{\Bbb C}$ is the parametrized real curve
$$
S^1 \longrightarrow \Bbb R^2,\quad \roman e^{i\alpha} \mapsto
(u(\alpha),v(\alpha))\in \Bbb R^2, \ (\alpha \in \Bbb R),
\tag3.10.3
$$
where $u(\alpha) + i v(\alpha) = 2 \roman e^{i\alpha} + \roman
e^{-2i\alpha}$; here $\Bbb R^2$ is embedded into $\Bbb C^2$ as the
real affine space of real points in $\Bbb C^2$ in the obvious
fashion. The curve (3.10.3) is a {\it hypocycloid\/}, as noted in
\cite\chakirus\ (Section 5). The real orbit space $\roman
{SU}(3)\big/\roman {SU}(3)\cong T\big/S_3 $ relative to conjugation, realized
within the model $T^{\Bbb C}\big/S_3$ for the categorical quotient
$\roman {SL}(3,\Bbb C)\big/\big/\roman {SL}(3,\Bbb C)$, amounts to
the compact region in $\Bbb R^2$ enclosed by the curve (3.10.3).

As a complex algebraic stratified K\"ahler space, the adjoint
quotient looks considerably more complicated. Indeed, the
complexification $\Bbb R[T^{\Bbb C}\big/S_3]_{\Bbb C}$ of the real
coordinate ring $\Bbb R[T^{\Bbb C}\big/S_3]\cong \Bbb R[T^{\Bbb
C}]^{S_3}$ of the adjoint quotient $T^{\Bbb C}\big/S_3$, viewed as
a real semialgebraic space, has the nine generators
(3.7.3), subject to the relations (3.7.5), (3.7.8), and (3.7.10),
together with the relations $\sigma_3 = 1$ and $\overline \sigma_3
= 1$. In view of the last two relations, it suffices to take the
generators
$$
\sigma_1,\ \overline \sigma_1, \ \sigma_2,\ \overline \sigma_2, \
\sigma, \ \rho, \ \overline \rho, \tag3.10.4
$$
subject to the relations which arise from the relations (3.7.5),
(3.7.8), and (3.7.10) by substitution of $1$ for each occurrence of
 $\sigma_3$ and $\overline \sigma_3$.
This does not change the relation (3.7.10) since neither $\sigma_3$
nor $\overline \sigma_3$ occur in (3.7.10) and, applied to the
relations (3.7.5) and (3.7.8), the procedure yields the relations
$$
\aligned
D_3(1,-\sigma_1,\sigma_2, - 1)
\overline \sigma_2
&=(9  + \sigma_2^3 - 4 \sigma_1 \sigma_2)
\overline \sigma_1^2
\\
&
\quad
+(4 \sigma_1^2  -3 \sigma_2 - \sigma_1 \sigma_2^2)
\overline \sigma_1 \sigma
\\
&
\quad
+(6 \sigma_1  - \sigma_2^2) \overline \sigma_1 \rho
+(\sigma_2^2-3\sigma_1  )  \sigma^2
\\
&
\quad
+(9  - \sigma_1 \sigma_2) \sigma \rho
+(\sigma_1^2 - 3 \sigma_2) \rho^2
\endaligned
\tag3.10.5
$$
and
$$
\aligned
D_3(1,-\sigma_1,\sigma_2, - 1)
&=
 \overline \sigma_1^3
-\sigma_2 \overline \sigma_1^2 \sigma
+
(\sigma_2^2 - 2 \sigma_1) \overline \sigma_1^2 \rho
+ \sigma_1 \overline \sigma_1 \sigma^2
\\
&\quad
-((\sigma_1^2 - 2 \sigma_2)\sigma_1
- \sigma_1^3 + 3 \sigma_1 \sigma_2 - 3)
\overline \sigma_1 \sigma \rho
\\
&\quad
-\sigma^3
+(\sigma_1^2 - 2 \sigma_2) \overline \sigma_1 \rho^2
\\
&\quad+ \sigma_2 \sigma^2 \rho
- \sigma_1 \sigma \rho^2
+ \rho^3 .
\endaligned
\tag3.10.6
$$
Thus the complexification $\Bbb R[T^{\Bbb C}\big/S_3]_{\Bbb C}$ of
the real coordinate ring of $T^{\Bbb C}\big/S_3$ is generated by
the seven bisymmetric functions (3.10.4), subject to the three
relations (3.7.10), (3.10.5) and (3.10.6).

Theorem 3.8.2 yields the stratified symplectic Poisson structure
on the complexification of the real coordinate ring of $T^{\Bbb
C}\big/S_3$ in terms of the nine  generators $\tau_{(1,0)}$,
$\tau_{(2,0)}$, $\tau_{(3,0)}$, $\tau_{(0,1)}$, $\tau_{(0,2)}$,
$\tau_{(0,3)}$, $\tau_{(1,1)}$, $\tau_{(1,2)}$, $\tau_{(2,1)}$.
This Poisson structure has rank 4 on the top stratum, rank 2 on
the middle stratum, and rank zero at the three points of the
bottom stratum, that is, at the three cusps of the complex affine
curve (3.10.2).

Finally, cf. (3.6) above, the reduced K\"ahler potential
$\kappa_{\roman{red}}$ is the function in $C^{\omega}(T^{\Bbb
C}\big/S_3) \cong C^{\omega}(T^{\Bbb C})^{S_3}$ which arises from
the K\"ahler potential $\kappa$, cf. (1.2) and (2.1.3), when
$\kappa$ is rewritten as a real analytic function of (the real and
imaginary parts of) the generators (3.10.4), similarly as in
(3.9.6) above. We refrain from spelling out the details.

\medskip\noindent {\bf 4. Other simple compact Lie groups}
\smallskip\noindent
We illustrate the theory developed so far briefly for various
other groups.

\smallskip\noindent
(4.1) {\smc  $B_n$: $K=\roman{SO}(2n+1, \Bbb R)$ ($n \geq 1$)\/}.
The adjoint quotient is isomorphic to that for the following case;
indeed, the root systems $B_n$ and $C_n$ differs from each other
only by an exchange of long and short roots.
The statement of the theorem in the introduction for
the group $\roman{SO}(2n+1,\Bbb R)$ is a consequence of
Corollary 3.4.4, combined with the fact that, under
the embedding of
this group into $\roman U(2n+1)$
given by the defining representation,
the first $n$ fundamental representations
of $\roman U(2n+1)$ restrict to the fundamental representations
of $\roman{SO}(2n+1,\Bbb R)$.

\smallskip\noindent
(4.2) {\smc  $C_n$: $K=\roman{Sp}(n)$ ($n \geq 1$)\/}. Embed the product
$\roman{Sp}(1)\times \dots\times \roman{Sp}(1)$ of $n$ copies
of $\roman{Sp}(1)$ into $\roman{Sp}(n)$ in the standard fashion.
Accordingly, the product of the standard maximal tori in these copies
of $\roman{Sp}(1)$ yields the standard maximal torus $T$ of $\roman{Sp}(n)$,
and the Weyl group $W$ amounts to the wreath product
$(\Bbb Z/2)^n \rtimes S_n$ of $\Bbb Z/2$
with the symmetric group $S_n$ on $n$ letters. Here each copy of $Z/2$
arises as the Weyl group of a copy of $\roman{Sp}(1)$ and,
with a grain of salt, the symmetric
group  $S_n$ yields the overall symmetries of the situation.
In terms of  the notation introduced in (2.2) above, we take
$V=\Bbb C^{2n}$ and embed $T^{\Bbb C}$ into $V$ in the obvious way, so that
the following hold:

\noindent
(i) The obvious coordinate functions
$z_1,\dots, z_n, z_1^{-1},\dots, z_n^{-1}$
on $T^{\Bbb C}$ generating the complex algebraic coordinate ring
$$
\Bbb C[T^{\Bbb C}]=\Bbb C[z_1,\dots, z_n, z_1^{-1},\dots,
z_n^{-1}] \tag4.2.1
$$ are independent coordinate functions on $V$.

\noindent
(ii) The action of the $j$'th copy of $\Bbb Z/2$ in $W$ interchanges
$z_j$ and $z_j^{-1}$ and leaves invariant the other coordinates,
and the copy of $S_n$ in $W$ permutes the
$z_1,\dots, z_n$'s and the $z_1^{-1},\dots, z_n^{-1}$'s
separately; this specifies the action of the Weyl group $W$.

Extending the notation introduced in (3.9) above, for $1 \leq j \leq n$,
let $Z_j = z_j + z_j^{-1}$. The algebra $\Bbb C[T^{\Bbb C}]^{(\Bbb Z/2)^n}$
of $(\Bbb Z/2)^n$-invariants is plainly the polynomial algebra
$\Bbb C[Z_1,\dots,Z_n]$, and the algebra
$\Bbb C[T^{\Bbb C}]^W$ of $W$-invariants is given by
$$
\Bbb C[T^{\Bbb C}]^W = \Bbb C[Z_1,\dots,Z_n]^{S_n}. \tag4.2.2
$$
This algebra of invariants is freely generated by the $n$
elementary symmetric functions or by the first $n$ power sum
functions in the variables $Z_1,\dots,Z_n$, and the resulting
Hilbert map from $T^{\Bbb C}$ to $\Bbb C^{n}$ induces a complex
algebraic isomorphism from the adjoint quotient $T^{\Bbb C}\big/W$
onto a copy of $\Bbb C^{n}$. This yields the complex algebraic
structure of the adjoint quotient $T^{\Bbb C}\big/W$. The $n$
elementary symmetric  functions in the variables $Z_1,\dots,Z_n$
yield the fundamental characters of $\roman{Sp}(n,\Bbb C)$.

To derive the real structure, we note first that
the complexification
$\Bbb R[T^{\Bbb C}]_{\Bbb C}$ of the real coordinate ring
$\Bbb R[T^{\Bbb C}]$ of $T^{\Bbb C}$ amounts to the complex algebra
$$
\Bbb R[T^{\Bbb C}]_{\Bbb C} = \Bbb C[z_1,\dots, z_n,
z_1^{-1},\dots, z_n^{-1}, \overline z_1,\dots, \overline z_n,
\overline z_1^{-1},\dots, \overline z_n^{-1}]. \tag4.2.3
$$
Extending the notation introduced in (3.9) above further, for $1
\leq j \leq n$, let
$$
\overline Z_j = \overline z_j + \overline z_j^{-1}, \ \sigma_{(j)}
= \frac{z_j}{\overline z_j} + \frac{\overline z_j}{z_j},\ \RRR_j =
\sigma_{(j)}^2 + Z_j^2 + \overline Z_j^2 - Z_j \overline Z_j
\sigma_{(j)} -4.
$$
In view of the discussion in (3.9) above,
cf. the defining relation (3.9.1), the algebra
of $(\Bbb Z/2)^n$-invariants is given by
$$
\Bbb R[T^{\Bbb C}]_{\Bbb C}^{(\Bbb Z/2)^n} = \Bbb
C[Z_1,\dots,Z_n,\overline Z_1,\dots,\overline Z_n,
\sigma_{(1)},\dots,\sigma_{(n)}]\big/(\RRR_1, \dots \RRR_n).
\tag4.2.4
$$
Hence the algebra $\Bbb R[T^{\Bbb C}]_{\Bbb C}^{W}$ of
$W$-invariants is the algebra of trisymmetric functions in the
variables
$$
Z_1,\dots,Z_n,\overline Z_1,\dots,\overline Z_n,
\sigma_{(1)},\dots,\sigma_{(n)},
$$
subject to the relations $\RRR_1=0$, \dots, $\RRR_n=0$, where the
symmetric group $S_n$ permutes the $Z_1,\dots,Z_n$'s, the
$\overline Z_1,\dots,\overline Z_n$'s, and the
$\sigma_{(1)},\dots,\sigma_{(n)}$'s separately in the obvious
fashion. Moreover, cf. (3.9.4), the Poisson structure is induced
by the identities
$$
\frac i2 \{Z_j, \overline Z_j\}= Z_j\overline Z_j -
2\sigma_{(j)},\ (1 \leq j \leq n). \tag4.2.5
$$
The statement of the theorem in the introduction for
the group $K=\roman{Sp}(n)$ is a consequence of
Corollary 3.4.4, combined with the fact that, under
the embedding of
this group into $\roman U(2n)$
given by the defining representation,
the first $n$ fundamental representations
of $\roman U(2n)$ restrict to the fundamental representations
of $\roman{Sp}(n)$.

\smallskip\noindent
(4.3) {\smc  $D_n$: $K=\roman{SO}(2n, \Bbb R)$ ($n \geq 2$)\/}.
The obvious embedding of the product of $n$ copies of
$\roman{SO}(2,\Bbb R)$ into $\roman{SO}(2n,\Bbb R)$ yields the
standard maximal torus $T$ of $\roman{SO}(2n,\Bbb R)$, and the
group $(\Bbb Z/2)^n \rtimes S_n$ acts on $T$ in the same fashion
as on the torus in (4.2) above. However, the Weyl group $W\cong
(\Bbb Z/2)^{n-1} \rtimes S_n$ amounts to  the subgroup of $(\Bbb
Z/2)^n \rtimes S_n$ where among the substitutions in the copy of
$(\Bbb Z/2)^n$ only even ones are admitted.

In terms of  the notation introduced in (2.2) above,
similarly as in (4.2) above, we take
$V=\Bbb C^{2n}$ and embed $T^{\Bbb C}$ into $V$ in the standard way, so that
the obvious coordinate functions
$z_1,\dots, z_n, z_1^{-1},\dots, z_n^{-1}$
on $T^{\Bbb C}$ generating the complex algebraic coordinate ring
$$
\Bbb C[T^{\Bbb C}]=\Bbb C[z_1,\dots, z_n, z_1^{-1},\dots, z_n^{-1}]
\tag4.3.1
$$
are independent coordinate functions on $V$ and so that what
corresponds to (4.2)(ii) still holds. As above, extending the
notation introduced in (3.9), for $1 \leq j \leq n$, let $Z_j =
z_j + z_j^{-1}$ and let $\sigma_j$ be the $j$'th elementary
symmetric functions  in the variables $Z_1,\dots,Z_n$. The $n$'th
elementary symmetric function $\sigma_n(Z_1,\ldots, Z_n)$
decomposes as a sum of $2^n$ terms of the kind $\frac{z_1 \ldots
z_k}{z_{k+1} \ldots z_n}$, and we refer to the difference $k-(n-k)
= 2k-n$ as the {\it degree\/} of $\frac{z_1 \ldots z_k}{z_{k+1}
\ldots z_n}$. Now $\sigma_n(Z_1,\ldots, Z_n)=Z_1 \cdot\ldots \cdot
Z_n$ decomposes as a sum
$$
Z_1 \cdot\ldots \cdot Z_n = \sigma^+_n(z_1,\dots, z_n,
z_1^{-1},\dots, z_n^{-1}) + \sigma^-_n(z_1,\dots, z_n,
z_1^{-1},\dots, z_n^{-1})
$$
where $\sigma^+_n$ is the sum of all terms of degrees congruent to
$n$ modulo 4 and $\sigma^-_n$
 the sum of all terms of degrees congruent to $n-2$ modulo 4.
For example, when $n=3$,
$$
\align \sigma^+_n&= z_1 z_2 z_3 + \frac {z_1}{z_2 z_3} +\frac
{z_2}{z_3 z_1} + \frac {z_3}{z_1 z_2},
\\
\sigma^-_n&= \frac{z_1 z_2}{z_3}  + \frac {z_2 z_3}{z_1} +\frac
{z_3 z_1}{z_2} + \frac 1{z_1 z_2 z_3}
\endalign
$$
The functions $\sigma^+_n$ and $\sigma^-_n$ are both $(\Bbb
Z/2)^{n-1}$-invariant, even $W$-invariant. The $n$'th exterior
power $\Lambda^n$ of the defining representation of
$\roman{SO}(2n+1,\Bbb R)$, restricted to $\roman{SO}(2n,\Bbb R)$
(relative to the obvious embedding), is well known to decompose as
a direct sum of two irreducible $\roman{SO}(2n,\Bbb
R)$-representations $\Lambda^+$ and $\Lambda^-$, and $\sigma^+_n$
and $\sigma^-_n$ are the  characters of these representations. The
algebra $\Bbb C[T^{\Bbb C}]^{(\Bbb Z/2)^{n-1}}$ of $(\Bbb
Z/2)^{n-1}$-invariants is the subalgebra of (4.3.1) generated by
$Z_1,\dots,Z_n$ and $\sigma^+_n$ (or $\sigma^-_n$), and these are
subject to the relation
$$
(\sigma_n^+ + 2 \sigma_{n-2} + \dots)
(\sigma_n^- + 2 \sigma_{n-2} + \dots)
= (\sigma_{n-1} + \dots )^2.
\tag4.3.2
$$
Consequently the complex coordinate ring $\Bbb C[T^{\Bbb C}\big/W]$
of the adjoint quotient $T^{\Bbb C}\big/W$ for
$\roman{SO}(2n,\Bbb R)$, that is,
the algebra $\Bbb C[T^{\Bbb C}]^W$ of $W$-invariants, is the
subalgebra of (4.3.1) generated by the
 $n-1$ elementary symmetric
functions $\sigma_1$, \dots, $\sigma_{n-1}$ together with
$\sigma^+_n$ and $\sigma^-_n$, subject to the relation (4.3.2).
That this relation is indeed defining reflects the standard
structure of the complex representation ring of
$\roman{SO}(2n,\Bbb R)$; it is also a consequence of the formulas
(4.5.1) and (4.5.2) for $\delta^2_+$, $\delta^2_-$, and $\delta_+
\delta_-$ in (4.5) below. For $1 \leq j \leq n-1$, when $\chi_j$,
refers to the $j$'th fundamental character of $\roman{SO}(2n,\Bbb
C)$, that is, to the character of the $j$'th exterior power of the
definig representation,  the function $\sigma_j$ coincides with
the character $\chi_{j-1}+ \chi_j$.

The real structure and the Poisson structure can then be determined
in a way similar to that explained in (4.2) above,
but with an appropriate algebra of multisymmetric functions.
This is a bit messy, but there is no real difficulty.
We spare the reader and ourselves these added troubles here.

\smallskip\noindent
(4.4) {\smc  $B_n$: $K=\roman{Spin}(2n+1, \Bbb R)$ ($n \geq
2$)\/}. The obvious embedding of the product of $n$ copies of
$\roman{SO}(2,\Bbb R)$ into $\roman{SO}(2n+1,\Bbb R)$ yields the
standard maximal torus $T$ of $\roman{SO}(2n+1,\Bbb R)$, and the
Weyl group $W\cong (\Bbb Z/2)^n \rtimes S_n$ acts on the maximal
torus $T$ of $\roman{SO}(2n+1,\Bbb R)$ in the same fashion as on
the torus in (4.2) above. We take as maximal torus $\widetilde T$
in $\roman{Spin}(2n+1, \Bbb R)$ the pre-image of $T$ under the
canonical surjection from $\roman{Spin}(2n+1, \Bbb R)$ to
$\roman{SO}(2n+1, \Bbb R)$, and we realize $\widetilde T$ as the
subspace
$$
\widetilde T = \{(z_1,\dots,z_n,z); z_1 \ldots z_n = z^2\}
\subseteq T \times S^1
$$
of $T \times S^1$. The action of the Weyl group $W$ on $\widetilde
T$, restricted to the symmetric group $S_n$, viewed as a subgroup
of $W$ in the obvious way, amounts simply to permutation of the
$z_1,\dots,z_n$'s while $z$ remains fixed and, for $1\leq j \leq
n$, the unique lift to $\widetilde T$ of the involution which on
$T$ sends the coordinate $z_j$ to $z_j^{-1}$ is now given by the
assignment to $(z_j,z)$ of $(z_j^{-1},zz_j^{-1})$ and leaves
invariant the other coordinates. The action of the group (of order
2)  of deck transformations is given by the assignment to $z$ of
$-z$ and leaves the coordinates $z_1$, \dots, $z_n$ invariant.

In terms of  the notation introduced in (2.2) above, similarly as
in (4.2) above, we take $V=\Bbb C^{2n+2}$ and embed $\widetilde
T^{\Bbb C}$ into $V$ in such a way that the obvious coordinate
functions $z_1,\dots, z_n, z, z_1^{-1},\dots, z_n^{-1},z^{-1} $ on
$\widetilde T^{\Bbb C}$ generating the complex algebraic
coordinate ring
$$
\Bbb C[\widetilde T^{\Bbb C}]=\Bbb C[z_1,\dots, z_n, z,
z_1^{-1},\dots, z_n^{-1},z^{-1}]\big/(z_1\dots z_nz^{-2} - 1,
z^{-1}_1\dots z^{-1}_nz^{2} - 1) \tag4.4.1
$$
are independent coordinate functions on $V$ and so that what
corresponds to (4.2)(ii) still holds. As above, extending the
notation introduced in (3.9), for $1 \leq j \leq n$, let $Z_j =
z_j + z_j^{-1}$, and let $\sigma_1,\dots,\sigma_n$ be the
elementary symmetric functions in the variables $Z_1,\dots,Z_n$.
Moreover, let
$$
\delta = z + \frac z{z_1} + \frac z{z_2} + \dots + z^{-1}
\tag4.4.2
$$
be the $(\Bbb Z/2)^n$-orbit sum of $z$; since $z$ is invariant
under $S_n$, $\delta$ is invariant even under $W$. A calculation
yields
$$
\delta^2 = 2^n + 2^{n-1} \sigma_1 + \dots + 2 \sigma_{n-1} +
\sigma_n. \tag4.4.3
$$
The algebra $\Bbb C[\widetilde T^{\Bbb C}]^{(\Bbb Z/2)^n}$ of $(\Bbb
Z/2)^n$-invariants is generated by $Z_1$, \dots, $Z_n$, and
$\delta$, subject to the relation (4.4.2), and the algebra $\Bbb
C[\widetilde T^{\Bbb C}]^W$ of $W$-invariants is the polynomial
algebra
$$
\Bbb C[\widetilde T^{\Bbb C}]^W = \Bbb
C[\sigma_1,\dots,\sigma_{n-1},\delta].
$$
We note that $\delta$ is the character of the half-spin
representation $\Delta$ of $\roman{Spin}(2n+1, \Bbb R)$ of
dimension $2^n$, and the identity (4.4.3) may be rewritten as the
familiar identity
$$
\Delta^2 = 1 + \Lambda^1 + \dots + \Lambda^n
\tag4.4.4
$$
in the complex representation ring of $\roman{Spin}(2n+1, \Bbb R)$ where
as usual $\Lambda^j$ ($1 \leq j \leq n)$ refers to the $j$'th
exterior power of the ordinary vector representation of dimension
$2n+1$ given by the projection to $\roman{SO}(2n+1, \Bbb R)$;
indeed, $\Lambda^1$ has character $\sigma_1 +1$, $\Lambda^2$ has
character $\sigma_2 +\sigma_1 +2$, etc.
 The subalgebra generated by
$\sigma_1,\dots,\sigma_n$ is precisely the isomorphic image in
$\Bbb C[\widetilde T^{\Bbb C}\big/W] \cong\Bbb C[\widetilde
T^{\Bbb C}]^W$ of the complex coordinate ring $\Bbb C[T^{\Bbb
C}\big/W]\cong \Bbb C[T^{\Bbb C}]^W$ of the adjoint quotient
$T^{\Bbb C}\big/W$ for $\roman{SO}(2n+1, \Bbb R)$ under the
canonical injection induced by the covering projection from
$\widetilde T$ to $T$, the algebra $\Bbb C[T^{\Bbb C}]^W$ having
been spelled out explicitly in (4.2) above. This algebra plainly
coincides with the subalgebra $\Bbb C[\widetilde T^{\Bbb
C}\big/W]^{\Bbb Z/2}$ of invariants under the induced action of
the group ($\cong \Bbb Z/2$) of deck transformations.

A variant of the approach in (4.2) above yields the real
semialgebraic structure and the stratified K\"ahler structure: The
complexification $\Bbb R[\widetilde T^{\Bbb C}]_{\Bbb C}$ of the
real coordinate ring $\Bbb R[\widetilde T^{\Bbb C}]$ of
$\widetilde T^{\Bbb C}$ amounts to the complex algebra
$$
\Bbb R[\widetilde T^{\Bbb C}]_{\Bbb C} = \Bbb C[z_1,\dots, z_n,
z_1^{-1},\dots, z_n^{-1}, \overline z_1,\dots, \overline z_n,
\overline z_1^{-1},\dots, \overline
z_n^{-1},z,z^{-1}]\big/(S,T,\overline S,\overline T),
$$
where $S=z_1\dots z_nz^{-2} - 1$ and  $T=z^{-1}_1\dots
z^{-1}_nz^{2} - 1$, cf. (4.4.1). Let $\sigma_{(n+1)}$ be the
$(\Bbb Z/2)^n$-orbit sum of $\frac z{\overline z}$. With the notation
introduced above, the algebra $\Bbb R[\widetilde T^{\Bbb C}]_{\Bbb
C}^{(\Bbb Z/2)^n}$ of $(\Bbb Z/2)^n$-invariants is generated by
$$
Z_1,\dots,Z_n,\delta,\overline Z_1,\dots,\overline
Z_n,\overline\delta, \sigma_{(1)},\ldots,\sigma_{(n)},
\sigma_{(n+1)}, \{\delta,\overline \delta\},
$$
subject to suitable relations, and the induced $S_n$-action
permutes the $Z_1,\dots,Z_n$'s, the $\overline Z_1,\dots,\overline
Z_n$'s, and the $\sigma_{(1)},\ldots,\sigma_{(n)}$'s separately
and leaves $\delta$, $\overline \delta$, $\sigma_{(n+1)}$,
and $\{\delta,\overline \delta\}$ invariant.
We shall comment on the Poisson bracket
$\{\delta,\overline \delta\}$ below.
 Hence the algebra $\Bbb R[\widetilde T^{\Bbb C}]_{\Bbb
C}^{W}$ of $W$-invariants arises from the algebra of trisymmetric
functions in the variables
$$
Z_1,\dots,Z_n,\overline Z_1,\dots,\overline Z_n,
\sigma_{(1)},\dots,\sigma_{(n)},
\tag4.4.5
$$
together with the four invariants $\delta$, $\overline \delta$,
$\{\delta, \overline \delta\}$ and
$\sigma_{(n+1)}$.

We will now explain how the stratified Poisson structure
can be determined.
The Poisson structure on the complexification
$\Bbb R[\widetilde T^{\Bbb C}]_{\Bbb C}$
of the real coordinate ring of $\widetilde T^{\Bbb C}$
is given by the formulas
$$
\align
\frac i2 \{z_j,\overline z_j\}&= 2 z_j \overline z_j\ (1 \leq j \leq n),
\tag4.4.6
\\
\frac i2 \{z,\overline z\}&=  z \overline z,
\tag4.4.7
\\
\frac i2 \{z_j,\overline z\}&= z_j \overline z \ (1 \leq j \leq n),
\tag4.4.8
\endalign
$$
similar to the formulas (4.2.5).
The factor 2 in (4.4.6) has been introduced for convenience,
in particular to arrive at simple formulas;
without the factor 2 in (4.4.6) we would need the factor
$\frac 12$ in (4.4.7) and (4.4.8).
This fixes of course the K\"ahler potential $\kappa$ on
$K^{\Bbb C}$.
A tedious but straightforward calculation yields
uniquely determined $W$-invariants $A_n$ and $B_n$ such that
$$
\frac i2\{\delta,\overline \delta\} = A_n - 2 \sigma_{(n+1)}
\tag4.4.9
$$
and
$$
\delta \overline \delta = A_n + B_n.
\tag4.4.10
$$
For example, when $n=1$ so that
$\roman{Spin}(3,\Bbb R) \cong \roman{SU}(2)$,
$B_1=0$,
$$
\delta \overline \delta = A_1 = (z+z^{-1})(\overline z + \overline
z^{-1}),
$$
and (4.4.9) comes down to
(3.9.4).
Likewise, when $n=2$,
$$
A_2 = (z+z^{-1})(\overline z + \overline
z^{-1})+
\left(\frac z{z_1}+\frac z{z_2}\right)
\left(\frac {\overline z}{\overline z_1}+\frac {\overline z}{\overline
z_2}\right).
$$
It is now straightforward to
 complete the determination of the stratified K\"ahler
structure on the adjoint quotient under discussion. We do not
pursue this here.

This discussion reveals that
the statement of the theorem in the introduction
is not true
for the group $K= \roman{Spin}(2n+1,\Bbb R)$ when $n \geq 2$.
Indeed, as a module over  the complexification of the real coordinate
ring of the adjoint quotient for $\roman{SO}(2n+1,\Bbb R)$,
the algebra of $W$-invariants
$\Bbb R[\widetilde T^{\Bbb C}]_{\Bbb
C}^{W}$ is generated by $\delta$, $\overline\delta$,
$\{\delta,\overline \delta\}$
and $\sigma_{(n+1)}$,
but {\it not\/}
by $\delta$, $\overline\delta$,
$\{\delta,\overline \delta\}$
alone.

\smallskip\noindent
(4.5) {\smc  $D_n$: $K=\roman{Spin}(2n, \Bbb R)$ ($n \geq 2$)\/}.
Relative to the standard embedding of $\roman{Spin}(2n, \Bbb R)$
into $\roman{Spin}(2n+1, \Bbb R)$, the maximal torus $\widetilde
T$ in $\roman{Spin}(2n+1, \Bbb R)$ is a maximal torus in
$\roman{Spin}(2n, \Bbb R)$. As in (4.3) above, the Weyl group
$W\cong (\Bbb Z/2)^{n-1} \rtimes S_n$ amounts to the subgroup of
$(\Bbb Z/2)^n \rtimes S_n$ where among the substitutions in the
copy of $(\Bbb Z/2)^n$ only even ones are admitted. The invariant
$\delta$ introduced in (4.4) above decomposes as a sum ${ \delta =
\delta^+ + \delta^- }$ where $\delta^+$ and $\delta^-$ are the
characters of the positive and negative half-spin representations
$\Delta^+$ and $\Delta^-$, respectively, of $\roman{Spin}(2n, \Bbb
R)$ of dimension $2^{n-1}$, and the algebra $\Bbb C[\widetilde
T^{\Bbb C}]^W$ of $W$-invariants is the polynomial algebra
$$
\Bbb C[\widetilde T^{\Bbb C}]^W = \Bbb
C[\sigma_1,\dots,\sigma_{n-2},\delta^+, \delta^-].
$$
Moreover, calculations yield the familiar identities in the
representation ring of $\roman{Spin}(2n, \Bbb R)$
 relating the
products $\delta^+ \delta^+$, $\delta^+ \delta^-$,
$\delta^-\delta^-$ with the elementary symmetric functions
$\sigma_1,\dots,\sigma_{n-1}$ and the functions $\sigma_n^+$ and
$\sigma_n^-$.
For $n$ odd,
$$
\aligned
\delta_+^2 &= \sigma_+ + 2 \sigma_{n-2} + \dots + 2^{n-2}\sigma_1
\\
\delta_-^2 &= \sigma_- + 2 \sigma_{n-2} + \dots + 2^{n-2}\sigma_1
\\
\delta_+ \delta_- &=
\sigma_{n-1} +2 \sigma_{n-3} + \dots + 2^{n-1}
\endaligned
\tag4.5.1
$$
while for $n$ even,
$$
\aligned
\delta_+^2 &= \sigma_+ + 2 \sigma_{n-2} + \dots + 2^{n-1}
\\
\delta_-^2 &= \sigma_- + 2 \sigma_{n-2} + \dots + 2^{n-1}
\\
\delta_+ \delta_- &=
\sigma_{n-1} +2 \sigma_{n-3} + \dots + 2^{n-2}\sigma_1
\endaligned
\tag4.5.2
$$

The non-trivial involution in the group ($\cong \Bbb Z/2$) of deck
transformations sends $\delta^+$ and $\delta^-$ to $-\delta^+$ and
$ -\delta^-$, respectively, and the subalgebra $\Bbb C[\widetilde
T^{\Bbb C}\big/W]^{\Bbb Z/2}$ of invariants under the induced
action of this group is the subalgebra generated by
$\sigma_1,\dots,\sigma_{n-1},\sigma_n^+,\sigma_n^-$. This algebra
is precisely the algebra $\Bbb C[T^{\Bbb C}]^W$ spelled out in
(4.3) above; it manifestly coincides with the isomorphic image in
$\Bbb C[\widetilde T^{\Bbb C}\big/W] \cong\Bbb C[\widetilde
T^{\Bbb C}]^W$ of the complex coordinate ring $\Bbb C[T^{\Bbb
C}\big/W]\cong \Bbb C[T^{\Bbb C}]^W$ of the adjoint quotient
$T^{\Bbb C}\big/W$ for $\roman{SO}(2n, \Bbb R)$ under the
canonical injection induced by the covering projection from
$\widetilde T$ to $T$.

The  stratified K\"ahler structure can then be determined
explicitly in a way similar to that hinted at in (4.4) above. We
do not pursue this here. We also note that the statement of the theorem
in the introduction is not true for $K=\roman{Spin}(2n,\Bbb R)$
when $n \geq 4$.

\smallskip\noindent
(4.6) {\smc  $K=G_{2(-14)}$\/}. This is the group of automorphisms
of the octonions $\Bbb O$. The defining representation has real
dimension 7; it is the subspace $\Bbb O_0$ of trace zero
octonions. The real dimension of $G_{2(-14)}$ equals 14. The long
roots constitute the root system $A_2$ and, accordingly,
$G_{2(-14)}$ contains a copy of $\roman{SU}(3)$ which, in turn,
admits an obvious interpretation in terms of the geometry of the
octonions. The Weyl group $W$ being a dihedral group ($\cong \Bbb
Z/6 \rtimes \Bbb Z/2$) of order 12, decomposes as $W\cong
W_{\roman{long}} \rtimes W_{\roman{short}}$, where
$W_{\roman{long}}\cong S_3$ is the Weyl group of the subsystem of
long roots and where $W_{\roman{short}}\cong \Bbb Z/2$  may be
taken as being generated by the reflection in the short simple
root. The standard maximal torus
$$
T =\{\roman{diag}(z_1,z_2,z_3); \, z_1 z_2 z_3 = 1\}
$$
for $\roman{SU}(3)$ is a maximal torus of $G_{2(-14)}$ as well
and, in terms of this torus, the Weyl group is generated by the
Weyl group $W_{\roman{long}}\cong S_3$ of $\roman{SU}(3)$ together
with the involution which sends $(z_1,z_2,z_3)$ to
$(z^{-1}_1,z^{-1}_2,z^{-1}_3)$. Indeed, this involution realizes
the obvious outer automorphism of $\roman{SU}(3)$ which
corresponds to the obvious symmetry of the root system $A_2$
interchanging the two roots. With the notation
$\sigma_1=z_1+z_2+z_3$ and $\sigma_2 = z_1 z_2 + z_1 z_3 + z_2
z_3$ used earlier, the complex coordinate ring $\Bbb C[T^{\Bbb
C}/S_3]$ of the adjoint quotient for $\roman{SU}(3)$ is the
complex polynomial algebra generated by $\sigma_1$ and $\sigma_2$,
and these are the characters of the defining representation (say)
$V_3\cong \Bbb C^3$ of $\roman{SL}(3,\Bbb C)$ and of the second
exterior square $\Lambda^2_{\Bbb C}V_3\cong \overline V_3$ of the
defining representation, respectively. Indeed, the outer
automorphism of $\roman{SU}(3)$ corresponding to the symmetry of
the root system interchanges the two fundamental representations
and passes to the involution on the adjoint quotient $T^{\Bbb
C}\big/S_3\cong \Bbb C^2$ for $\roman{SU}(3)$ which interchanges
the complex coordinate functions $\sigma_1$ and $\sigma_2$. Hence
the complex coordinate ring  of the adjoint quotient $T^{\Bbb
C}\big/W$ for $G_{2(-14)}$ is the algebra
$$
\Bbb C[T^{\Bbb C}/W] = \Bbb C[\sigma_1,\sigma_2]^{\Bbb Z/2} = \Bbb
C[\Sigma_1,\Sigma_2]
$$
where $\Sigma_1=\sigma_1 +\sigma_2$ and $\Sigma_2 =\sigma_1
\sigma_2$. In particular, complex algebraically, the adjoint
quotient $T^{\Bbb C}\big/W$ for $G_{2(-14)}$ is again a copy of
$\Bbb C^2$.

To interpret this description in terms of the two fundamental
characters of $G_{2(-14)}$, we note that, as a real
$\roman{SU}(3)$-representation, the defining representation of
$G_{2(-14)}$ decomposes as $V_3\oplus \Bbb R$ where $V_3$ is
viewed as a real 6-dimensional representation of $\roman{SU}(3)$.
Hence the defining representation has as character the function
$\Sigma_1 +1$. The other fundamental representation of
$G_{2(-14)}$ is the adjoint representation. Since the tensor
product $V_3 \otimes \Lambda^2V_3 \cong V_3 \otimes \overline V_3$
of the two fundamental representations of $\roman{SU}(3)$
decomposes as the direct sum $\frak {su}(3) \oplus \Bbb R$ of the
adjoint representation and the trivial 1-dimensional
representation, the adjoint representation of $\roman{SU}(3)$  has
character $\sigma_1 \sigma_2-1$. As a real
$\roman{SU}(3)$-representation, the adjoint representation of
$G_{2(-14)}$ decomposes as $\frak{su}(3) \oplus V_3$. Hence the
adjoint representation of $G_{2(-14)}$ has as character the
function $\Sigma_1 +\Sigma_2-1$. We note that, as a complex
$\roman{SU}(3)$-representation, the adjoint
 representation of $G_{2(-14)}$ decomposes as $
\frak {sl}(3,\Bbb C) \oplus V_3 \oplus \overline V_3 $.

The real semialgebraic structure and the stratified K\"ahler
structure arise from that spelled out in (3.10) above by taking
$(\Bbb Z/2)$-invariants. In view of the discussion in (3.10)
above, the complexification $\Bbb R[T^{\Bbb
C}\big/W_{\roman{long}}]_{\Bbb C}$ of the real coordinate ring of
the adjoint quotient $T^{\Bbb C}\big/W_{\roman{long}}$ for
$\roman{SU}(3)$ is generated by the seven invariants $\sigma_1$,
$\overline \sigma_1$, $\sigma_2$, $\overline \sigma_2$, $\sigma$,
$\rho$,  $\overline \rho$, cf. (3.10.4), subject to the three
relations (3.7.10), (3.10.5) and (3.10.6). The non-trivial
involution coming from the induced action of $W_{\roman{short}}
\cong \Bbb Z/2$ on $\Bbb R[T^{\Bbb C}\big/W_{\roman{long}}]_{\Bbb
C}$ interchanges, respectively, $\sigma_1$, $\overline\sigma_1$,
$\rho$ and $\sigma_2$, $\overline \sigma_2$, $\overline \rho$, and
sends $\sigma$ to a $W_{\roman{long}}$-invariant function (say)
$\overline \sigma$ on $T^{\Bbb C}$ which, on the real torus $T$,
viewed as a subspace of $T^{\Bbb C}$, coincides with $\sigma$. The
real semialgebraic structure and the Poisson structure can then be
determined in a way similar to that explained in (4.2) above, but
with an appropriate algebra of quatrisymmetric functions, relative
to the group $W_{\roman{short}} \cong \Bbb Z/2$.

Similarly as before, the statement of the theorem in the
introduction for the group $K=G_{2(-14)}$ is a consequence of
Corollary 3.4.4, combined with the observation that, under the
embedding of $G_{2(-14)}$ into $\roman U(7)$ given by the defining
representation $\Bbb O_0$ of $G_{2(-14)}$, (i) {\sl the first
fundamental\/} (i.~e. defining) {\sl representation of $\roman
U(7)$ restricts to the first fundamental\/} (i.~e. defining) {\sl
representation of $G_{2(-14)}$\/}, and (ii) {\sl the second
fundamental representation of\/} $\roman U(7)$ (the second
exterior square of the defining representation of $\roman U(7)$)
{\sl restricts to the sum of the two fundamental complex
representations $\frak g_{2}(\Bbb C)$ and $\Bbb O_0\otimes \Bbb C$
of $G_{2(-14)}$\/}. Indeed, (i) is obvious. To justify (ii), we
note first that the defining representation $\Bbb O_0$ of
$G_{2(-14)}$ is real and embeds $G_{2(-14)}$ into
$\roman{SO}(7,\Bbb R)$ and, as a $G_{2(-14)}$-representation, the
adjoint representation (Lie algebra) $\frak{so}(7,\Bbb R)$
decomposes as
$$
\frak{so}(7,\Bbb R) \cong \frak g_{2(-14)} \oplus \Bbb O_0.
$$
The adjoint representation of $\roman{SO}(7,\Bbb R)$ is the second
exterior square of the defining representation of
$\roman{SO}(7,\Bbb R)$ whence (ii).

\smallskip\noindent
(4.7) {\smc  $K=F_{4(-52)}$\/}. This is the group of automorphisms
of the exceptional Jordan algebra $\Cal H_3(\Bbb O)$ of hermitian
$(3 \times 3)$-matrices over the octonions $\Bbb O$. The defining
representation has real dimension 26; it is the subspace of trace
zero matrices in $\Cal H_3(\Bbb O)$. The defining representation
may also be obtained from a reductive decomposition of $\frak
e_{6(-78)}$ of the kind $\frak e_{6(-78)}= \frak f_4 + \frak p$.
The real dimension of $F_{4(-52)}$ equals 52. The long roots
constitute the root system $D_4$ and, accordingly, $F_{4(-52)}$
contains a copy of $\roman{Spin}(8,\Bbb R)$ which, in turn, admits
an obvious interpretation in terms of the geometry of the
octonions:  Denote by $V_8$ the (real) vector representation and
by $S_+$ and $S_-$ the (real) positive and negative spinor
representations, respectively, of $\roman{Spin}(8,\Bbb R)$. The
division algebra structure on $\Bbb O$ corresponds to a {\it
triality\/} $V_8\oplus S_+ \oplus S_- \to \Bbb R$, and
$\roman{Spin}(8,\Bbb R)$ is the symmetry group of this triality,
cf. \cite\adabotwo. The subgroup of $\roman{Spin}(8,\Bbb R)$ which
consists of the symmetries of the triality which are automorphism
of $\Bbb O$ is precisely the group $G_{2(-14)}$. The Weyl group of
$W$ of $F_{4(-52)}$ decomposes as $((\Bbb Z/3)^3\rtimes \Bbb S_4)
\rtimes S_3$ where the copy of $(\Bbb Z/3)^3\rtimes \Bbb S_4$ is
the Weyl group $W_{\roman{long}}$ of the system of long roots,
that is, of $\roman{Spin}(8,\Bbb R)$.
 The  maximal torus
$$
\widetilde T \cong\{(z_1,z_2,z_3,z_4,z); \, z_1 z_2 z_3 z_4= z^2\}
\subseteq (S^1)^5
$$
for $\roman{Spin}(8,\Bbb R)$ spelled out in (4.5) above for the
general case of $\roman{Spin}(2n,\Bbb R)$
 is a maximal
torus of $F_{4(-52)}$ as well and, in terms of this torus, the
Weyl group is generated by the Weyl group  of $\roman{Spin}(8,\Bbb
R)$ together with certain involutions; these involutions come from
the obvious outer automorphisms of $\roman{Spin}(8,\Bbb R)$ which
correspond to the obvious symmetry group $S_3$  of the root system
$D_4$ permuting the three boundary roots. This symmetry group was
discovered already by E. Cartan. It is generated by the
reflections in the short simple roots, and we will denote it by
$W_{\roman{short}}$. The vector representation  and the positive
and negative spinor representations  of $\roman{Spin}(8,\Bbb R)$
are all isomorphic under these outer automorphisms of
$\roman{Spin}(8,\Bbb R)$.

 With the notation $\sigma_1$, $\sigma_2$ for
the first two elementary symmetric functions in the variables
$Z_1,Z_2,Z_3,Z_4$ and the notation $\delta^+$, and $\delta^-$
 introduced in (4.5) above,
 the complex coordinate ring $\Bbb C[T^{\Bbb
C}/W_{\roman{long}}]$ of the adjoint quotient for
$\roman{Spin}(8,\Bbb R)$ is the complex polynomial algebra
generated by $\sigma_1$, $\sigma_2$, $\delta^+$, and $\delta^-$.
Here $\delta^+$ and $\delta^-$ are the characters of the positive
and negative spinor representations and, when $\chi_1$ refers to
the character of the vector representation $V_8$ of
$\roman{Spin}(8,\Bbb R)$  and $\chi_2$ to that of the adjoint
representation,  $\sigma_1$ coincides with $\chi_1+1$ and
$\sigma_2$ with $\chi_1+ \chi_2$.

The $W_{\roman{short}}$-action passes to an $S_3$-action on the
adjoint quotient $T^{\Bbb C}\big/W_{\roman{long}}\cong \Bbb C^4$
for $\roman{Spin}(8,\Bbb R)$ which permutes the complex coordinate
functions $\chi_1$, $\delta^+$, and $\delta^-$,
 whence the complex coordinate ring of the
adjoint quotient $T^{\Bbb C}\big/W$ for $F_{4(-52)}$ is the
algebra
$$
\Bbb C[T^{\Bbb C}/W] = \Bbb C[\chi_1,\chi_2, \delta^+, \delta^-
]^{S_3} = \Bbb C[\Sigma_1, \chi_2,\Sigma_2, \Sigma_3],
$$
where $\Sigma_1$, $\Sigma_2$, and $\Sigma_3$ are the elementary
symmetric functions in the variables $\chi_1, \delta^+, \delta^-$.
Since, as a $\roman{Spin}(8,\Bbb R)$-representation,
$$
\Cal H_3(\Bbb O)\cong \Bbb R^3 \oplus V_8 \oplus
 S_+ \oplus S_-,
$$
the character of the defining representation of $F_{4(-52)}$ is
the function $\Sigma_1 +2$. Furthermore, since, as a
$\roman{Spin}(8,\Bbb R)$-representation,
$$
\frak f_{4(-52)} \cong \frak {so}(8,\Bbb R)\oplus V_8 \oplus S_+
\oplus S_-,
$$
the sum
$$
\sigma_2 + \delta = \Sigma_1 + \chi_2
$$
is the character of the adjoint representation of $F_{4(-52)}$.
This is the second fundamental representation of $F_{4(-52)}$. We
do not make explicit the relationship between the virtual
characters $ \Sigma_1, \chi_2,\Sigma_2, \Sigma_3$ and the two
remaining fundamental representations of $F_{4(-52)}$.

In particular, complex algebraically, the adjoint quotient
$T^{\Bbb C}\big/W$ for $F_{4(-52)}$ is the quotient of
the adjoint quotient $T^{\Bbb
C}\big/W_{\roman{long}}\cong \Bbb C^4$ for $\roman{Spin}(8,\Bbb
R)$ relative to the induced $S_3$-action, and this quotient is
again a copy of $\Bbb C^4$.
Moreover,
the real semialgebraic structure arises from that hinted at in
(4.5) in the special case where $n=4$ by the operation of taking
$S_3$-invariants and, likewise,
the stratified K\"ahler structure can be determined explicitly by
an application of  the procedure of taking $S_3$-invariants to the
stratified K\"ahler structure on the adjoint quotient $T^{\Bbb
C}\big/W_{\roman{long}}\cong \Bbb C^4$ for $\roman{Spin}(8,\Bbb
R)$ which, in turn, was hinted at in (4.5) above. We do not pursue this
here.

\beginsection 5. Energy quantization on the adjoint quotient

Choose a {\it dominant\/} Weyl chamber in the Lie algebra $\frak
t$ of the maximal torus $T$ of $K$ and let $R^+$ be the resulting
system of positive roots. Let $\Delta_K$ denote the {\it
Casimir\/} operator on $K$ associated with the bi-invariant
Riemannian metric on $K$, and let $m=\dim K$. When $X_1,\dots,X_m$
is an orthonormal basis of $\frak k$,
$$
\Delta_K = X^2_1 +\dots + X^2_m
$$
in the universal algebra $\roman U(\frak k)$ of $\frak k$, cf.
e.~g. \cite\enelson\ (p.~591). The Casimir operator depends only
on the Riemannian metric, though. Since the metric on $K$ is
bi-invariant, so is the operator $\Delta_K$; hence, by Schur's
lemma, for each highest weight $\lambda$, the isotypical $(K\times
K)$-summand $L^2(K,dx)_\lambda \subseteq L^2(K,dx)$ associated
with $\lambda$ in the Peter-Weyl decomposition of $L^2(K,dx)$ is
an eigenspace, whence the representative functions are
eigenfunctions for $\Delta_K$. Let $\rho = \frac 12 \sum_{\alpha
\in R^+} \alpha$, so that $2\rho$ is the sum of the positive
roots. The eigenvalue of $\Delta_K$ corresponding to the highest
weight $\lambda$ is known to be given explicitly by $-\varepsilon$
where
$$\varepsilon_{\lambda}=(|\lambda+\rho|^2-|\rho|^2),
$$
cf. e.~g. \cite\helbotwo\ (Ch. V.1 (16) p.~502). The present sign
is dictated by the interpretation in terms of the energy spelled
out below. Thus $\Delta_K$ acts on each isotypical summand
$L^2(K,dx)_\lambda$ as scalar multiplication by
$-\varepsilon_{\lambda}$. The Casimir operator is known to
coincide with the nonpositive Laplace-Beltrami operator associated
with the (bi-invariant) Riemannian metric on $K$, see e.~g.
\cite\taylothr\ (A 1.2). In the Schr\"odinger picture (vertical
quantization on $\roman T^*K$), the operator $\widehat E_K$ which
arises as the unique extension of $-\frac 12 \Delta_K$ to an
unbounded self-adjoint operator on $L^2(K,dx)$ is the quantum
mechanical {\it energy\/} operator associated with the Riemannian
metric, whence the spectral decomposition of this operator refines
in the standard manner to the Peter-Weyl decomposition of
$L^2(K,dx)$ into isotypical $(K\times K)$-summands. The energy
operator may by obtained by vertical quantization of the geodesic
flow.

Let $\varepsilon$ be the symplectic volume form on $\roman
T^*K\cong K^{\Bbb C}$ (inducing Liouville measure). Define the
function $ \eta \colon K^{\Bbb C} \longrightarrow \Bbb R$ by
$$
\eta(x,Y) = \left(\roman{det}\left(\frac{\sin(\roman{ad}(Y))}
{\roman{ad}(Y)}\right)\right)^{\frac 12}, \ x \in K, \,Y \in \frak
k; \tag5.1
$$
this yields a non-negative real analytic function on $K^{\Bbb C}$
which depends only on the variable $Y \in \frak k$ and, for $x\in
K$ and $Y\in \frak k$, we will also write $\eta(Y)$ instead of
$\eta(x,Y)$. The function $\eta^2$ is the density of Haar measure
relative to the Liouville volume measure on $K^{\Bbb C}$, see
e.~g. \cite\bhallfou\ (Lemma 5). Both measures are
$K$-bi-invariant; in particular, as a function on $\frak k$,
$\eta$ is $\roman{Ad}(K)$-invariant.

Half-form K\"ahler quantization on $\roman T^*K\cong K^{\Bbb C}$
is accomplished by the Hilbert space
 $\Cal
H L^2(K^{\Bbb C},\roman e^{-\kappa} \ETAA \varepsilon)$ of
holomorphic functions which are square integrable relative to the
measure $\roman e^{-\kappa} \ETAA \varepsilon$ \cite\bhallone,
\cite\holopewe.
 Via the embedding of $\frak k$ into $\frak k^{\Bbb
C}$, the operator $\Delta_K$ is a differential operator on
$K^{\Bbb C}$. In view of Theorem 5.2 in \cite\holopewe\ which, in
turn, is a consequence of the holomorphic Peter-Weyl theorem
\cite\holopewe, in the holomorphic quantization on $\roman
T^*K\cong K^{\Bbb C}$, the unique extension $\widehat E_{K^{\Bbb
C}}$ of  $-\frac 12 \Delta_K$ to an (unbounded) self-adjoint
operator on $\Cal H L^2(K^{\Bbb C},\roman e^{-\kappa} \ETAA
\varepsilon)$ is the quantum mechanical {\it energy\/} operator
associated with the Riemannian metric, and the spectral
decomposition of this operator refines to the holomorphic
Peter-Weyl decomposition of $\Cal H L^2(K^{\Bbb C},\roman
e^{-\kappa} \ETAA \varepsilon)$ into isotypical $(K\times
K)$-summands; see \cite\holopewe\ (Sections 6 and 7) for details.

As before let  $T^{\Bbb C}\subseteq K^{\Bbb C}$ be the
complexification of $T$ and $W$  the Weyl group. We will now
explain how
 {\sl half-form K\"ahler quantization of the reduced kinetic energy
$\frac 12 \kappa_{\roman{red}}$ on the adjoint quotient of
$K^{\Bbb C}$ leads to the the irreducible characters of\/}
$K^{\Bbb C}$. To this end, let $dt$ be Haar measure on $T^{\Bbb
C}$, and let $d[t]$ be the measure on $T^{\Bbb C}/W$ which,
multiplied by the order $|W|$ of the Weyl group $W$, is the push
forward of the measure $dt$ on $T^{\Bbb C}$ under the projection
from $T^{\Bbb C}$ to $T^{\Bbb C}/W$. The Haar measure $dt$
actually coincides with Liouville measure on $T^{\Bbb C}$ and, on
the regular part of $T^{\Bbb C}\big/W$, viewed as an ordinary
smooth symplectic manifold, $d[t]$ coincides with Liouville
measure; we will therefore write $\varepsilon_{\roman{red}}$ for
$d[t]$ as well. Given a highest weight $\lambda$ for $K$, let
$\chi_{\lambda}\colon K^{\Bbb C} \to \Bbb C$ be the irreducible
(algebraic) character of $K^{\Bbb C}$ with highest weight
$\lambda$. Each such character $\chi_{\lambda}$ of $K^{\Bbb C}$
manifestly passes to an algebraic function on the adjoint quotient
$T^{\Bbb C}\big/ W$ and we denote this function by
$[\chi_{\lambda}]$.

The {\it saturation\/} of the zero locus $\mu^{-1}(0)$ is the
$K^{\Bbb C}$-closure $K^{\Bbb C}\mu^{-1}(0) \subseteq K^{\Bbb C}$
of $\mu^{-1}(0)$ in $K^{\Bbb C}$ relative to the (conjugation)
action of $K^{\Bbb C}$ on itself, and the inclusion
$\mu^{-1}(0)\subseteq K^{\Bbb C}\mu^{-1}(0)$ induces a
homeomorphism from the reduced space $(\roman T^*K)_0 =\mu^{-1}(0)
\big/ K$ onto the $K^{\Bbb C}$-quotient $K^{\Bbb C}\mu^{-1}(0)
\big /K^{\Bbb C}$. This yields an alternate description of the
quotient $K^{\Bbb C}\big/\big/K^{\Bbb C} \cong T^{\Bbb C}/W$ of
$K^{\Bbb C}$. Relative to the projection to the quotient, consider
the push forward to $T^{\Bbb C}/W$ of the measure $ \roman
e^{-\kappa} \eta \varepsilon =\frac {\roman e^{-\kappa}}{\eta} dx
dY$ on $K^{\Bbb C}$. This push forward measure on $T^{\Bbb C}/W$
has a density relative to Liouville measure
$\varepsilon_{\roman{red}}$ and hence can be written in the form
$$
\roman e^{-\kappa_{\roman{red}}}\gamma\varepsilon_{\roman{red}},
\tag5.2
$$
for a uniquely determined real valued function $\gamma$ on
$T^{\Bbb C}\big/W$ such that, given two holomorphic $K$-invariant
functions $\Phi$ and $\Psi$ on $K^{\Bbb C}$ that are square
integrable relative to the measure $ \roman e^{-\kappa} \eta
\varepsilon$, when $\Phi_{\roman{red}}$ and $\Psi_{\roman{red}}$
denote the induced holomorphic functions on the quotient $T^{\Bbb
C}\big/W$,
$$
\int_{K^{\Bbb C}} \overline \Phi \Psi \roman e^{-\kappa} \eta
\varepsilon =\int_{T^{\Bbb C}\big/W} \overline
{\Phi_{\roman{red}}} \Psi_{\roman{red}} \roman
e^{-\kappa_{\roman{red}}}\gamma\varepsilon_{\roman{red}}. \tag5.3
$$
To establish the existence of the function $\gamma$ consider the
conjugation mapping
$$
q^{\Bbb C}\colon \left(K^{\Bbb C}\big/ T^{\Bbb C}\right) \times
T^{\Bbb C} \longrightarrow K^{\Bbb C},\ (y T^{\Bbb C},t)\mapsto
yty^{-1},\  y\in K^{\Bbb C}, t\in T^{\Bbb C},
$$
and integrate the induced $(2m)$-form $(q^{\Bbb C})^*(\roman
e^{-\kappa} \eta \varepsilon)$ over \lq\lq the fiber\rq\rq\
$K^{\Bbb C}\big/ T^{\Bbb C}$. In view of the Gaussian constituent
$\roman e^{-\kappa}$, this integration is a well defined
operation. Let $\widetilde \gamma$ be the density of  the
resulting $(2n)$-form on $T^{\Bbb C}$  relative to the Liouville
volume form on $T^{\Bbb C}\cong \roman T^*T$ where $n=\dim T$ and
let $\widehat \gamma = \widetilde \gamma\big/|W|$ where $|W|$
denotes the order of the Weyl group $W$. The function $\widehat
\gamma$ descends to a function on the quotient $T^{\Bbb C}\big/W$,
and dividing this function by $\roman e^{-\kappa_{\roman{red}}}$
we obtain the function $\gamma$ we are looking for.

The half-form quantization procedure on $K^{\Bbb C}$ induces a
half-form quantization procedure on the adjoint quotient $T^{\Bbb
C}\big/W$; we do not spell out the details here.  The above
discussion then leads to the following result, which includes the
quantization of the reduced kinetic energy $\frac
12\kappa_{\roman{red}}\in C^\omega(T^{\Bbb C}\big/W)$.

\proclaim{Theorem 5.4} The  quantum Hilbert space for the
stratified K\"ahler structure on the adjoint quotient $T^{\Bbb
C}\big/W$ amounts to the Hilbert space ${ \Cal H L^2(T^{\Bbb
C}\big/W,\roman e^{-\kappa_{\roman{red}}} \gamma
\varepsilon_{\roman{red}}) }$ of holomorphic functions on $T^{\Bbb
C}\big/W$ that are square-integrable with respect to the measure
$\roman e^{-\kappa_{\roman{red}}} \gamma
\varepsilon_{\roman{red}}$, and this Hilbert space is freely
spanned by the holomorphic functions $[\chi_{\lambda}]$ on
$T^{\Bbb C}\big/W$ that correspond to the irreducible characters
of $K^{\Bbb C}$. Furthermore,  the reduced energy operator
$\widehat E_{\roman{red}}$ is given by
$$
\widehat E_{\roman{red}}[\chi_{\lambda}] =
\varepsilon_{\lambda}[\chi_{\lambda}].
$$
\endproclaim

Indeed, let $\varepsilon_T$ be the Liouville volume form on
$T^{\Bbb C}\cong \roman T^*T$. The restriction mapping induces a
unitary isomorphism
$$
\Cal H L^2(K^{\Bbb C},\roman e^{-\kappa} \eta \varepsilon)^K
\longrightarrow
 \Cal H L^2(T^{\Bbb C}, {\roman e^{-\kappa}
\widehat\gamma \varepsilon_T})^W
$$
of Hilbert spaces, and the canonical map
$$
\Cal H L^2(T^{\Bbb C}\big/W,\roman e^{-\kappa_{\roman{red}}}
\gamma \varepsilon_{\roman{red}}) \longrightarrow \Cal H
L^2(T^{\Bbb C}, {\roman e^{-\kappa} \widehat\gamma
\varepsilon_T})^W
$$
is an isomorphism of Hilbert spaces. We shall give elsewhere an
intrinsic description of the quantization downstairs on $\Cal H
L^2(T^{\Bbb C}\big/W,\roman e^{-\kappa_{\roman{red}}} \gamma
\varepsilon_{\roman{red}})$. The advantage of the description in
terms of the quotient $T^{\Bbb C}\big/W$ rather than in terms of
the $W$-invariants of $\Cal H L^2(T^{\Bbb C},\roman e^{-\kappa}
\widetilde\gamma \varepsilon_T)$ is that it brings the {\it
costratified\/} nature of the Hilbert space structure to the fore:

The holomorphic quantization procedure in \cite\qr\ yields a {\it
costratified Hilbert space\/}, that is, a system of Hilbert
spaces, one Hilbert space for the closure of each stratum, with
bounded linear operators among these Hilbert spaces which
correspond to the closure relations among the strata. Such a
system is {\sl structure on the quantum level which has the
classical singularities as its shadow\/}.
 Under the present
circumstances, the costratified structure arises in the following
fashion: As explained before, the adjoint quotient $T^{\Bbb
C}\big/W$ is decomposed into strata, and the closure of a stratum
is an affine complex variety. Given a stratum $S$, let $I_S$ be
the ideal of functions in the complex coordinate ring $\Bbb
C[T^{\Bbb C}\big/W]$ which vanish on $S$ or, equivalently, on its
closure $\overline S$, let $\widehat I_S \subseteq  \Cal H
L^2(T^{\Bbb C}\big/W,\roman e^{-\kappa_{\roman{red}}} \gamma
\varepsilon_{\roman{red}})$ be the closed subspace spanned by
$I_S$, and let $\Cal H_S$ be the orthogonal complement of
$\widehat I_S$ in $\Cal H L^2(T^{\Bbb C}\big/W,\roman
e^{-\kappa_{\roman{red}}} \gamma \varepsilon_{\roman{red}})$.
Equivalently, we can characterize $\widehat I_S$ as the space of
holomorphic functions in $\Cal H L^2(T^{\Bbb C}\big/W,\roman
e^{-\kappa_{\roman{red}}} \gamma \varepsilon_{\roman{red}})$ that
vanish on $S$ and, as a complex vector space, $\Cal H_S$ is the
space of holomorphic functions on the complex manifold $S$ which
arise as restrictions of holomorphic functions on $T^{\Bbb
C}\big/W$ which are square-integrable relative to the measure
$\roman e^{-\kappa_{\roman{red}}} \gamma
\varepsilon_{\roman{red}}$. In particular, when $S$ is the top
stratum, $\Cal H_S$ coincides with the entire Hilbert space $\Cal
H L^2(T^{\Bbb C}\big/W,\roman e^{-\kappa_{\roman{red}}} \gamma
\varepsilon_{\roman{red}})$.
 Thus, for each stratum $S$,
there is a canonical projection from $\Cal H L^2(T^{\Bbb
C}\big/W,\roman e^{-\kappa_{\roman{red}}} \gamma
\varepsilon_{\roman{red}})$ to $\Cal H_S$ and, given two strata
$S_1$ and $S_2$ with $S_2 \subseteq \overline {S_1}$, there is a
unique projection $\Cal H_{S_1}\to \Cal H_{S_2}$. The system
$\{\Cal H_S\}$, as $S$ ranges over the strata, together with the
projections $\Cal H_{S_1}\to \Cal H_{S_2}$ whenever $S_2 \subseteq
\overline {S_1}$, constitutes a {\it costratified\/} Hilbert
space. Alternatively, the quantum Hilbert space is spanned by the
functions $[\chi_{\lambda}]\in\Bbb C[T^{\Bbb C}\big/W]$ induced by
the irreducible characters as $\lambda$ ranges over the highest
weights. The closure of a stratum is a complex affine variety in
$T^{\Bbb C}\big/W$ ($\cong \Bbb C^{n}$ when $K$ is semisimple and
simply connected of rank $n$), and the Hilbert space corresponding
to that stratum is simply the quotient of $\Cal H L^2(T^{\Bbb
C}\big/W, \roman e^{- \kappa_{\roman{red}}}
\gamma\varepsilon_{\roman{red}})$ obtained when the functions
$[\chi_{\lambda}]$ are restricted to that stratum. In
\cite\hurusch, the quantum mechanics on such a costratified
Hilbert space is worked out in detail for the special case where
$K=\roman{SU}(2)$.

Thus the measure  giving rise to the Hilbert space
 $\Cal H L^2(T^{\Bbb C}\big/W, \roman e^{-
\kappa_{\roman{red}}} \gamma\varepsilon_{\roman{red}})$ on the
reduced level, viz. the measure $\roman e^{- \kappa_{\roman{red}}}
\gamma\varepsilon_{\roman{red}}$, appears to be formally of the
same kind as the measure $\roman e^{- \kappa} \eta\varepsilon$
which determines the unreduced Hilbert space $\Cal H L^2(K^{\Bbb
C}, \roman e^{- \kappa} \eta\varepsilon)$, that is to say, on the
reduced level, the measure involves the reduced energy
$\kappa_{\roman{red}}$, the reduced Liouville volume form
$\varepsilon_{\roman{red}}$, and a certain correction term
$\gamma$. However, there is a fundamental difference: On the
unreduced level, the correction term $\eta$ comes from the
metaplectic correction and is intrinsically defined in terms of
the geometry of the group $K^{\Bbb C}$ whereas the correction term
$\gamma$ on the reduced level is not merely defined in terms of
the geometry of the adjoint quotient $K^{\Bbb C}\big/\big/ K^{\Bbb
C} \cong T^{\Bbb C}\big/W$ and, in a sense, encapsulates part of
the history as to how the quotient arises. Relative to these
measures, in the case at hand, {\sl half-form K\"ahler
quantization commutes with reduction, observables and Hilbert
space structures included\/}, the reduced space being endowed with
a measure which is not in an obvious way related with merely the
geometry of the reduced space.

\medskip

\centerline{\smc References} \smallskip
\widestnumber\key{999}

\ref \no \adabotwo \by J. F. Adams \book Lectures on Exceptional
Lie groups \bookinfo eds. Z. Mahmud and M. Mimura \publ University
of Chicago Press \publaddr Chicago \yr 1996
\endref

\ref \no \armcusgo \by J. M. Arms,  R. Cushman, and M. J. Gotay
\paper  A universal reduction procedure for Hamiltonian group
actions \paperinfo in: The geometry of Hamiltonian systems, T.
Ratiu, ed. \jour MSRI Publ. \vol 20 \pages 33--51 \yr 1991 \publ
Springer Verlag \publaddr Berlin $\cdot$ Heidelberg $\cdot$ New
York $\cdot$ Tokyo
\endref

\ref \no \rbieltwo \by R. Bielawski \paper K\"ahler metrics on
$G^{\Bbb C}$ \jour J. reine angew. Mathematik \vol 559 \yr 2003
\pages 123--136 \finalinfo{\tt math.DG/0202255}
\endref

\ref \no \biemitwo \by E. Bierstone and P. D. Milman \paper
Semianalytic and subanalytic sets \jour Publ. Math. I. H. E. S.
\vol 67 \yr 1988 \pages  5--42
\endref

\ref \no \broetomd \by  Th. Br\"ocker and T. tom Dieck \book
Representations of Compact Lie groups \bookinfo Graduate Texts in
Mathematics, No. 98 \publ Springer Verlag \publaddr Berlin $\cdot$
Heidelberg  $\cdot$   New York $\cdot$  Tokyo \yr 1985
\endref

\ref \no \bourbalg
\by  N. Bourbaki
\book Algebra
\publ Springer Verlag
\publaddr Berlin  $\cdot$  Heidelberg  $\cdot$   New York $\cdot$  Tokyo
\yr 1989
\endref

\ref \no \chakirus \by S. Charzy\'nski, J. Kijowski, G. Rudolph,
and M. Schmidt \paper On the stratified classical configuration
space of lattice QCD \jour J. Geom. and Physics \vol 55 \yr 2005
\pages 137--178
\endref

\ref \no \cushbate \by R. H. Cushman and L. M. Bates \book Global
aspects of classical integrable systems \publ Birkh\"auser Verlag
\publaddr Boston $\cdot$ Basel $\cdot$ Berlin \yr 1997
\endref

\ref \no \dalbec \by J. Dalbec \paper Multisymmetric functions
\jour Beitr\"age zur Algebra und Geometrie, Contributions to Algebra
and Geometry \vol 40 \yr 1999 \pages 27--51
\endref

\ref\no\diracone \by P. A. M. Dirac \book Lectures on Quantum
Mechanics \publ Belfer Graduate School of Science,Yeshiva
University \publaddr New York \yr 1964
\endref

\ref \no \elhalnee \by  R. Elmore, P. Hall, and A. Neeman \paper
An application of classical invariant theory to identifiability in
nonparametric mixtures \jour Ann. Inst. Fourier \vol 55 \yr 2005
\pages 1--28
\endref

\ref \no \bhallfou \by  B. C. Hall \paper Phase space bounds for
quantum mechanics on a compact Lie group
 \jour Comm. in Math. Physics \vol 184 \yr 1997 \pages
233--250
\endref

\ref \no \bhallone \by  B. C. Hall \paper Geometric quantization
and the generalized Segal-Bargmann transform for Lie groups of
compact type \jour Comm. in Math. Physics \vol 226 \yr 2002 \pages
233--268 \finalinfo {\tt quant.ph/0012015}
\endref

\ref \no \helbotwo \by S. Helgason \book Groups and Geometric
Analysis. Integral geometry, invariant differential operators, and
spherical functions \bookinfo Pure and Applied Mathematics, vol.
113 \publ Academic Press Inc. \publaddr Orlando, Fl. \yr 1984
\endref

\ref \no \poiscoho \by J. Huebschmann \paper Poisson cohomology
and quantization \jour J. reine angew.
 Mathematik
\vol  408 \yr 1990 \pages 57--113
\endref

\ref \no  \souriau \by J. Huebschmann \paper On the quantization
of Poisson algebras \paperinfo Symplectic Geometry and
Mathematical Physics, Actes du colloque en l'honneur de Jean-Marie
Souriau, P. Donato, C. Duval, J. Elhadad, G.M. Tuynman, eds.;
Progress in Mathematics, Vol. 99 \publ Birkh\"auser Verlag
\publaddr Boston $\cdot$ Basel $\cdot$ Berlin \yr 1991 \pages
204--233
\endref

\ref \no \kaehler \by J. Huebschmann \paper K\"ahler spaces,
nilpotent orbits, and singular reduction \linebreak
\jour Memoirs
AMS \vol 172/814 \yr 2004 \publ Amer. Math. Society \publaddr
Providence, R. I.\finalinfo {\tt math.DG/0104213}
\endref

\ref \no \lradq \by J. Huebschmann \paper Lie-Rinehart algebras,
descent, and quantization \paperinfo Galois theory, Hopf algebras,
and semiabelian categories \jour Fields Institute Communications
\vol 43 \publ Amer. Math. Society \publaddr Providence, R.~I.
\pages 295--316 \yr 2004 \finalinfo {\tt math.SG/0303016}
\endref

\ref \no \qr \by J. Huebschmann \paper K\"ahler quantization and
reduction \finalinfo {\tt math.SG/0207166} \jour J.  reine angew.
Mathematik
\vol 591
\yr 2006
\pages 75--109
\endref

\ref \no \varna \by J. Huebschmann \paper Classical phase space
singularities and quantization \paperinfo in: Quantum Theory and
Symmetries. IV. V. Dobrev, ed. (to appear) \publ Heron Press
\publaddr Sofia \yr 2006 \finalinfo {\tt math-ph/0610047}
\endref

\ref \no \bedlepro \by J. Huebschmann \paper Singular
Poisson-K\"ahler geometry of certain adjoint quotients\paperinfo
in: The Mathematical Legacy of C. Ehresmann, J. Kubarski, and R.
Wolak, eds. \jour Banach Center Publications (to appear)
\finalinfo {\tt math.SG/0610614}
\endref

\ref \no \holopewe \by J. Huebschmann \paper The holomorphic
Peter-Weyl theorem and the Blattner-Kostant-Sternberg pairing
\finalinfo{\tt math.DG/0610613}
\endref

\ref \no \hurusch \by J. Huebschmann, G. Rudolph, and M. Schmidt
\paper A lattice gauge model for singular quantum mechanics
\paperinfo in preparation
\endref

\ref \no \humphone \by J. E. Humphreys \book Conjugacy classes in
semisimple algebraic groups \bookinfo Mathematical Surveys and
Monographs, vol. 43 \publ Amer. Math. Society \publaddr
Providence, R. I. \yr 1995
\endref

\ref \no \kemneone \by G. Kempf and L. Ness \paper The length of
vectors in representation spaces \jour Lecture Notes in
Mathematics \vol 732 \yr 1978 \pages 233--244 \paperinfo Algebraic
geometry, Copenhagen, 1978 \publ Springer Verlag \publaddr Berlin
$\cdot$ Heidelberg $\cdot$ New York
\endref

\ref \no \kirwaboo \by F. Kirwan \book Cohomology of quotients in
symplectic and algebraic geometry \publ Princeton University Press
\publaddr Princeton, New Jersey \yr 1984
\endref

\ref \no \lempszoe \by L. Lempert and R. Sz\"oke \paper Global
solutions of the homogeneous complex Monge-Am\-p\`ere equations
and complex structures on the tangent bundle of Riemannian
manifolds \jour Math. Ann. \vol 290 \yr 1991 \pages 689--712
\endref

\ref \no \lermonsj \by E. Lerman, R. Montgomery, and R. Sjamaar
\paper Examples of singular reduction \paperinfo Symplectic
Geometry, Warwick, 1990,  D. A. Salamon, editor, London Math. Soc.
Lecture Note Series, vol. 192 \yr 1993 \pages  127--155 \publ
Cambridge University Press \publaddr Cambridge, UK
\endref

\ref \no \lunatwo
\by D. Luna
\paper Sur certaines op\'erations diff\'erentiables des groupes de Lie
\jour Amer. J. of Math.
\vol 97
\yr 1975
\pages 172--181
\endref

\ref \no \lunathr
\by D. Luna
\paper Fonctions diff\'erentiables invariantes sous l'op\'eration
d'un groupe de Lie r\'eductif
\jour Ann. Inst. Fourier
\vol 26
\yr 1976
\pages 33--49
\endref

\ref \no \enelson \by E. Nelson \paper Analytic vectors \jour Ann.
of Mathematics \vol 70 \yr 1959 \pages  572--615
\endref

\ref \no \netto \by E. Netto \book Vorlesungen \"uber Algebra
\publ Teubner Verlag \publaddr Leipzig \yr 1896
\endref

\ref \no \procschw \by C. Procesi and G. Schwarz \paper
Inequalities defining orbit spaces \jour Invent. math. \vol 81 \yr
1985 \pages 539--554
\endref

\ref \no \richaone \by R. W. Richardson Jr. \paper Conjugacy
classes of $n$-tuples in Lie algebras and algebraic groups \jour
Duke Math. J. \vol 57 \yr 1988 \pages 1--35
\endref

\ref \no \gwschwar \by G. W. Schwarz \paper Smooth functions
invariant under the action of a compact Lie group \jour Topology
\vol 14 \yr 1975 \pages 63--68
\endref

\ref \no \gwschwat
\by G. W. Schwarz \paper The topology of
algebraic quotients \paperinfo In: Topological methods in
algebraic transformation groups, Progress in Mathematics, Vol. 80
\yr 1989 \pages 135--152 \publ Birkh\"auser Verlag\publaddr Boston
$\cdot$ Basel $\cdot$ Berlin
\endref

\ref \no \slodoboo \by P. Slodowy \book Simple singularities and
simple algebraic groups \bookinfo Lecture Notes in Mathematics,
Vol. 815 \yr 1980 \publ Springer \publaddr Berlin $\cdot$
Heidelberg $\cdot$ New York
\endref

\ref \no \snowone \by D. M. Snow \paper Reductive group actions on
Stein spaces \jour Math. Ann. \vol 259 \yr 1982 \pages 79--97
\endref

\ref \no \steinbtw \by R. Steinberg \paper Regular elements of
semisimple algebraic groups \jour Pub. Math. I. H. E. S. \vol 25
\yr 1965 \pages 49--80
\endref

\ref \no \szoekone \by R. Sz\"oke \paper Complex structures on
tangent bundles of Riemannian manifolds \jour Math. Ann. \vol 291
\yr 1991 \pages 409--428
\endref

\ref \no \taylothr \by J.  Taylor \paper The Iwasawa decomposition
and limiting behaviour of Brownian motion on symmetric spaces of
non-compact type \jour Cont. Math. \vol 73 \yr 1988 \pages
303--331
\endref

\ref \no \vaccarin \by F. Vaccarino \paper The ring of
multisymmetric functions \jour Ann. Inst. Fourier \vol 55 \yr 2005
\pages 717--731
\endref

\ref \no \weylbook \by H. Weyl \book The classical groups \publ
Princeton University  Press \publaddr Princeton, New Jersey \yr
1946
\endref

\enddocument